# Fourier series (based) multiscale method for computational analysis in science and engineering:

# II. Composite Fourier series method for simultaneous approximation to functions and their (partial) derivatives

Weiming Sun[a,*,+] and Zimao Zhang[b]

**Abstract:** Fourier series multiscale method, a concise and efficient analytical approach for multiscale computation, will be developed out of this series of papers. The second paper is concerned with simultaneous approximation to functions and their (partial) derivatives. On the basis of sufficient conditions of $2r$ ($r$ is a positive integer) times term-by-term differentiation of Fourier series, a one-dimensional or two-dimensional function with general boundary conditions is decomposed into the linear combination of a corner function that describes the discontinuities at corners of the domain (only for the two-dimensional function), boundary functions that describe the discontinuities on boundaries of the domain and an internal function that describes the smoothness within the domain. It leads to the methodology of simultaneous approximation of functions and their (partial) derivatives with composite Fourier series. And specifically, for the algebraical polynomial based composite Fourier series method, the reproducing property of complete algebraical polynomials is theoretically analyzed and the approximation accuracy is validated with numerical examples. This study generalizes the Fourier series method with supplementary terms to simultaneous approximation of functions and their (partial) derivatives up to $2r$th order.

**Keywords:** composite Fourier series; simultaneous approximation; (partial) derivatives

## 1. Introduction

Often in mathematics or engineering applications, we are faced with the task to approximate functions and their (partial) derivatives simultaneously by (trial) approximating functions and their corresponding (partial) derivatives. In some cases, algebraical polynomials [1-3] are employed for the approximating functions. In other more complicated cases, the approximating functions are constructed by various interpolation methods. These lead to the simultaneous approximation of functions and their (partial) derivatives by means of Lagrange interpolation polynomials [4-6], Hermite interpolation polynomials [7], Pal-type interpolation polynomials [8], or interpolation splines [9, 10]. However, when we solve boundary value problems governed by partial differential equations and boundary conditions, the approximating (solution) functions expanded in terms of Fourier series are preferred [11].

---

[a]Department of Mathematics and Big Data, School of Artificial Intelligence, Jianghan University, Wuhan, 430056, China
[b]Department of Mechanics, School of Civil Engineering, Beijing Jiaotong University, Beijing, 100044, China
*Correspondence to: Weiming Sun, Department of Mathematics and Big Data, School of Artificial Intelligence, Jianghan University, Wuhan, 430056, China
[+] E-mail: xuxinenglish@hust.edu.cn

For brevity, we denote the array of a function and its (partial) derivatives by the (partial) derivative vector of the function. Then for a specific (partial) derivative vector, the usual Fourier series approximation method, also referred to as the Fourier series direct-expansion method [12-23], is to directly expand the function, together with its various order (partial) derivatives, into their corresponding Fourier series to obtain a simultaneous approximation of the (partial) derivative vector. However, this approximation method is sensitive to the form of the function, and usually results in significant difference of convergence speeds in different parts (for example, the internal part and boundary part) of the domain and between different components of the (partial) derivative vector. All these properties revealed are internal unfavorable factors for the Fourier series direct-expansion method [24].

It's worthy of notice that formulating a function with a linear combination of a Fourier series and some supplementary terms presents a new approach to the simultaneous approximation of functions and their (partial) derivatives [25-59]. In this method, the supplementary terms in the linear combination are introduced to take care of potential discontinuities with the function and its (partial) derivatives on the boundaries when they are extended onto the entire $x_1$ axis (or $x_1 - x_2$ plane), and usually in simple forms such as algebraical polynomials for one-dimensional problems or single Fourier series for two-dimensional problems. Consequently, the Fourier series in the linear combination is only intended for the description of the smoothness behaviour of the function within the domain, and is uniformly convergent and termwise differentiable. It is a pity that the general formulas for Fourier series of higher order (partial) derivatives of the functions are heretofore unavailable, which restricts the application of the Fourier series method with supplementary terms to case studies concerning the simultaneous approximation of functions and their (partial) derivatives up to fourth order.

Therefore, in the second paper of the series of researches on Fourier series multiscale method, we take as a point of departure the formulas for Fourier series of higher order (partial) derivatives of one-dimensional and two-dimensional functions with general boundary conditions, and especially the sufficient conditions for $2r$ times ($r$ is a positive integer) term-by-term differentiation of Fourier series of functions [60], and dedicate ourselves to the improvement of the Fourier series method with supplementary terms. It is verified that the improved version of the Fourier series method with supplementary terms is feasible for functions with varied boundary conditions, and it also has excellent uniform and simultaneous convergence of functions and their (partial) derivatives up to $2r$th order. Furthermore, basis functions of different kinds in the approximation series, such as the functions in the trigonometrical system, the algebraical polynomials of specific degrees, and if necessary, the products of specific functions (usually the algebraical polynomials) and trigonometrical functions are in balanced use and indispensable. As this improved method develops an idea of composites of different kinds of basis functions, we rename it the composite Fourier series method. Accordingly, in this paper the structural decompositions of functions are described firstly. And then the framework of the composite Fourier series method is developed. Finally, detailed formulations related to the algebraical polynomial based composite Fourier series method are presented and validated with numerical examples.

## 2. Structural decompositions of functions

With the guidance of the sufficient conditions for $2r$ times ($r$ is a positive integer) term-by-term differentiation of Fourier series of functions, in this section we formulate a one-dimensional or two-dimensional function with general boundary conditions as a superposition with a corner function that describes the discontinuities at corners of the domain (only for the two-dimensional function), boundary functions that describe the discontinuities on boundaries of the domain and an internal function that describes the smoothness within the domain.

*2.1. For one-dimensional functions*

Given a positive integer $r$, we denote by $C^{2r-1}([-a,a])$ (or $C^{2r-1}([0,a])$) the set of $2r-1$ times continuously differentiable functions over the interval $[-a,a]$ (or $[0,a]$). Then the structural decomposition of the one-dimensional function $u(x_1)$ involves the following three situations:

1. $u(x_1) \in C^{2r-1}([-a,a])$, with its $2r$th order derivative absolutely integrable over the interval $[-a,a]$, and the resulting internal function is specified to satisfy the sufficient conditions for $2r$ times term-by-term differentiation of full-range Fourier series over the interval $[-a,a]$.

2. $u(x_1) \in C^{2r-1}([0,a])$, with its $2r$th order derivative absolutely integrable over the interval $[0,a]$, and the resulting internal function is specified to satisfy the sufficient conditions for $2r$ times term-by-term differentiation of half-range cosine series over the interval $[0,a]$.

3. $u(x_1) \in C^{2r-1}([0,a])$, with its $2r$th order derivative absolutely integrable over the interval $[0,a]$, and the resulting internal function is specified to satisfy the sufficient conditions for $2r$ times term-by-term differentiation of half-range sine series over the interval $[0,a]$.

For the first situation, we decompose the function $u(x_1)$ as follows:

1. Select the boundary function $\varphi_1(x_1)$ such that
$$\varphi_1^{(k_1)}(a) - \varphi_1^{(k_1)}(-a) = u^{(k_1)}(a) - u^{(k_1)}(-a), \quad k_1 = 0,1,\cdots,2r-1. \tag{1}$$

2. Let $\varphi_0(x_1) = u(x_1) - \varphi_1(x_1)$, then the internal function $\varphi_0(x_1)$ satisfies
$$\varphi_0^{(k_1)}(a) - \varphi_0^{(k_1)}(-a) = 0, \quad k_1 = 0,1,\cdots,2r-1. \tag{2}$$

Thus, we have decomposed the one-dimensional function $u(x_1)$ with general conditions into the sum of two functions, namely, the internal function $\varphi_0(x_1)$ and the boundary function $\varphi_1(x_1)$, where the internal function $\varphi_0(x_1)$ satisfies the sufficient conditions for $2r$ times term-by-term differentiation of the full-range Fourier series expansion of one-dimensional functions.

For the latter two situations, we perform the structural decompositions of the one-dimensional functions by similar procedures.

## 2.2. For two-dimensional functions

For a positive integer $r$, we use the notations $C^{2r-1}([-a,a]\times[-b,b])$ and $C^{2r-1}([0,a]\times[0,b])$ to denote the sets of $2r-1$ times continuously differentiable functions respectively over the domains $[-a,a]\times[-b,b]$ and $[0,a]\times[0,b]$. Then the structural decomposition of the two-dimensional function $u(x_1,x_2)$ involves the following five situations:

1. $u(x_1,x_2)$ is from $C^{2r-1}([-a,a]\times[-b,b])$ and all its $2r$ th order partial derivatives are absolutely integrable over the domain $[-a,a]\times[-b,b]$, and the resulting internal function is specified to satisfy the sufficient conditions for $2r$ times term-by-term differentiation of full-range Fourier series over the domain $[-a,a]\times[-b,b]$.

2. $u(x_1,x_2)$ is from $C^{2r-1}([0,a]\times[0,b])$ and all its $2r$ th order partial derivatives are absolutely integrable over the domain $[0,a]\times[0,b]$, and the resulting internal function is specified to satisfy the sufficient conditions for $2r$ times term-by-term differentiation of half-range cosine-cosine series over the domain $[0,a]\times[0,b]$.

3. $u(x_1,x_2)$ is from $C^{2r-1}([0,a]\times[0,b])$ and all its $2r$ th order partial derivatives are absolutely integrable over the domain $[0,a]\times[0,b]$, and the resulting internal function is specified to satisfy the sufficient conditions for $2r$ times term-by-term differentiation of half-range sine-cosine series over the domain $[0,a]\times[0,b]$.

4. $u(x_1,x_2)$ is from $C^{2r-1}([0,a]\times[0,b])$ and all its $2r$ th order partial derivatives are absolutely integrable over the domain $[0,a]\times[0,b]$, and the resulting internal function is specified to satisfy the sufficient conditions for $2r$ times term-by-term differentiation of half-range cosine-sine series over the domain $[0,a]\times[0,b]$.

5. $u(x_1,x_2)$ is from $C^{2r-1}([0,a]\times[0,b])$ and all its $2r$ th order partial derivatives are absolutely integrable over the domain $[0,a]\times[0,b]$, and the resulting internal function is specified to satisfy the sufficient conditions for $2r$ times term-by-term differentiation of half-range sine-sine series over the domain $[0,a]\times[0,b]$.

Here we give a detailed decomposition process as an illustration for the first situation.

1. Select the corner function $\varphi_3(x_1,x_2)$ such that
$$\varphi_3^{(k_1,k_2)}(a,b) - \varphi_3^{(k_1,k_2)}(a,-b) - \varphi_3^{(k_1,k_2)}(-a,b) + \varphi_3^{(k_1,k_2)}(-a,-b) =$$
$$u^{(k_1,k_2)}(a,b) - u^{(k_1,k_2)}(a,-b) - u^{(k_1,k_2)}(-a,b) + u^{(k_1,k_2)}(-a,-b),$$
$$k_1,k_2 = 0,1,2,\cdots,\quad k_1+k_2 \le 2r-2. \tag{3}$$

2. Let $\varphi_{012}(x_1,x_2) = u(x_1,x_2) - \varphi_3(x_1,x_2)$, then we have
$$\varphi_{012}^{(k_1,k_2)}(a,b) - \varphi_{012}^{(k_1,k_2)}(a,-b) - \varphi_{012}^{(k_1,k_2)}(-a,b) + \varphi_{012}^{(k_1,k_2)}(-a,-b) = 0,$$
$$k_1,k_2 = 0,1,2,\cdots,\quad k_1+k_2 \le 2r-2. \tag{4}$$

Further, we select the boundary function $\varphi_1(x_1,x_2)$ such that
$$\varphi_1^{(k_1,0)}(a,x_2) - \varphi_1^{(k_1,0)}(-a,x_2) = [u^{(k_1,0)}(a,x_2) - \varphi_3^{(k_1,0)}(a,x_2)]$$
$$-[u^{(k_1,0)}(-a,x_2) - \varphi_3^{(k_1,0)}(-a,x_2)],\quad x_2 \in (-b,b),\quad k_1 = 0,1,\cdots,2r-1, \tag{5}$$
$$\varphi_1^{(0,k_2)}(x_1,b) - \varphi_1^{(0,k_2)}(x_1,-b) = 0,\quad x_1 \in (-a,a),\quad k_2 = 0,1,\cdots,2r-1, \tag{6}$$
$$\varphi_1^{(k_1,k_2)}(a,b) - \varphi_1^{(k_1,k_2)}(a,-b) - \varphi_1^{(k_1,k_2)}(-a,b) + \varphi_1^{(k_1,k_2)}(-a,-b) = 0,$$
$$k_1,k_2 = 0,1,2,\cdots,\quad k_1+k_2 \le 2r-2. \tag{7}$$

Meanwhile, we select another boundary function $\varphi_2(x_1, x_2)$ such that

$$\varphi_2^{(k_1,0)}(a, x_2) - \varphi_2^{(k_1,0)}(-a, x_2) = 0, \quad x_2 \in (-b, b), \quad k_1 = 0, 1, \cdots, 2r-1, \tag{8}$$

$$\varphi_2^{(0,k_2)}(x_1, b) - \varphi_2^{(0,k_2)}(x_1, -b) = [u^{(0,k_2)}(x_1, b) - \varphi_3^{(0,k_2)}(x_1, b)]$$
$$-[u^{(0,k_2)}(x_1, -b) - \varphi_3^{(0,k_2)}(x_1, -b)], \quad x_1 \in (-a, a), \quad k_2 = 0, 1, \cdots, 2r-1, \tag{9}$$

$$\varphi_2^{(k_1,k_2)}(a, b) - \varphi_2^{(k_1,k_2)}(a, -b) - \varphi_2^{(k_1,k_2)}(-a, b) + \varphi_2^{(k_1,k_2)}(-a, -b) = 0,$$
$$k_1, k_2 = 0, 1, 2, \cdots, \quad k_1 + k_2 \leq 2r - 2. \tag{10}$$

3. Let the internal function $\varphi_0(x_1, x_2) = \varphi_{012}(x_1, x_2) - \varphi_1(x_1, x_2) - \varphi_2(x_1, x_2)$
$$= u(x_1, x_2) - \varphi_1(x_1, x_2) - \varphi_2(x_1, x_2) - \varphi_3(x_1, x_2), \tag{11}$$

then $\varphi_0(x_1, x_2)$ satisfies

$$\varphi_0^{(k_1,0)}(a, x_2) - \varphi_0^{(k_1,0)}(-a, x_2) = 0, \quad x_2 \in (-b, b), \quad k_1 = 0, 1, \cdots, 2r-1, \tag{12}$$

$$\varphi_0^{(0,k_2)}(x_1, b) - \varphi_0^{(0,k_2)}(x_1, -b) = 0, \quad x_1 \in (-a, a), \quad k_2 = 0, 1, \cdots, 2r-1, \tag{13}$$

$$\varphi_0^{(k_1,k_2)}(a, b) - \varphi_0^{(k_1,k_2)}(a, -b) - \varphi_0^{(k_1,k_2)}(-a, b) + \varphi_0^{(k_1,k_2)}(-a, -b) = 0,$$
$$k_1, k_2 = 0, 1, 2, \cdots, \quad k_1 + k_2 \leq 2r - 2. \tag{14}$$

Therefore, we can decompose the two-dimensional function $u(x_1, x_2)$ with general boundary conditions as the sum of the internal function $\varphi_0(x_1, x_2)$, the boundary functions $\varphi_1(x_1, x_2)$ and $\varphi_2(x_1, x_2)$ and the corner function $\varphi_3(x_1, x_2)$, where the internal function $\varphi_0(x_1, x_2)$ satisfies the sufficient conditions for $2r$ times term-by-term differentiation of full-range Fourier series of two-dimensional functions, the boundary function $\varphi_1(x_1, x_2)$ satisfies the sufficient conditions for $2r$ times term-by-term differentiation of full-range Fourier series of one-dimensional functions along the $x_2$-direction, and the boundary function $\varphi_2(x_1, x_2)$ satisfies the sufficient conditions for $2r$ times term-by-term differentiation of full-range Fourier series of one-dimensional functions along the $x_1$-direction.

## 3. Composite Fourier series method

On the basis of structural decompositions of functions described in the previous section, we select proper forms of corner functions (only for the two-dimensional functions), boundary functions and internal functions and develop the composite Fourier series method for simultaneous approximation of functions and their (partial) derivatives.

### 3.1. For one-dimensional functions

As an illustration, we derive the corresponding full-range composite Fourier series for the one-dimensional function $u(x_1)$ belonging to $C^{2r-1}([-a, a])$, along with its derivatives.

1. Constructing the boundary function $\varphi_1(x_1)$

We select the vector of supplementary functions as follows

$$\mathbf{p}_1^T(x_1) = [p_{1,1}(x_1) \quad p_{1,2}(x_1) \quad \cdots \quad p_{1,2r}(x_1)], \tag{15}$$

where $p_{1,j}(x_1)$, $j=1,2,\cdots,2r$, represent a set of closed-form sufficiently smooth functions defined over the interval $[-a,a]$. The term "sufficiently smooth" implies that $2r$th order derivatives of these functions exist and are continuous at any point over the interval $[-a,a]$. Such requirements can be readily satisfied by algebraical polynomials. Theoretically, there are an infinite number of choices for these supplementary functions.

Then for $k_1 = 0,1,\cdots,2r$, we define the vector of derivatives of the selected supplementary functions

$$\mathbf{p}_1^{(k_1)\text{T}}(x_1) = [p_{1,1}^{(k_1)}(x_1) \quad p_{1,2}^{(k_1)}(x_1) \quad \cdots \quad p_{1,2r}^{(k_1)}(x_1)]. \tag{16}$$

Further, let the vector of undetermined constants $\mathbf{a}_1^\text{T} = [G_{1,1} \quad G_{1,2} \quad \cdots \quad G_{1,2r}]$, then we can construct the boundary function as

$$\varphi_1(x_1) = \mathbf{p}_1^\text{T}(x_1) \cdot \mathbf{a}_1. \tag{17}$$

Considering the necessary conditions as given in Eq. (1) for the boundary function $\varphi_1(x_1)$, we have

$$\mathbf{R}_1 \mathbf{a}_1 = \mathbf{q}_1, \tag{18}$$

where

$$\mathbf{R}_1 = \begin{bmatrix} \mathbf{p}_1^{(0)\text{T}}(a) - \mathbf{p}_1^{(0)\text{T}}(-a) \\ \mathbf{p}_1^{(1)\text{T}}(a) - \mathbf{p}_1^{(1)\text{T}}(-a) \\ \vdots \\ \mathbf{p}_1^{(2r-1)\text{T}}(a) - \mathbf{p}_1^{(2r-1)\text{T}}(-a) \end{bmatrix}, \tag{19}$$

and the vector involving boundary values of the function $u(x_1)$ and its derivatives

$$\mathbf{q}_1^\text{T} = [u^{(0)}(a) - u^{(0)}(-a) \quad u^{(1)}(a) - u^{(1)}(-a) \quad \cdots \quad u^{(2r-1)}(a) - u^{(2r-1)}(-a)]. \tag{20}$$

Then we have

$$\mathbf{a}_1 = \mathbf{R}_1^{-1} \mathbf{q}_1. \tag{21}$$

If we denote the vector of basis functions by

$$\mathbf{\Phi}_1^\text{T}(x_1) = \mathbf{p}_1^\text{T}(x_1) \cdot \mathbf{R}_1^{-1}, \tag{22}$$

then the boundary function $\varphi_1(x_1)$ can finally be expressed as

$$\varphi_1(x_1) = \mathbf{\Phi}_1^\text{T}(x_1) \cdot \mathbf{q}_1. \tag{23}$$

2. Constructing the internal function $\varphi_0(x_1)$

Expand $\varphi_0(x_1)$ in a full-range Fourier series over the interval $[-a,a]$, and then we have

$$\varphi_0(x_1) = \sum_{m=0}^{\infty} \mu_m [V_{1m} \cos(\alpha_m x_1) + V_{2m} \sin(\alpha_m x_1)], \tag{24}$$

where $\alpha_m = m\pi/a$, $\mu_m = \begin{cases} 1/2, & m=0 \\ 1, & m>0 \end{cases}$, $V_{1m}$ and $V_{2m}$ are the Fourier coefficients of $\varphi_0(x_1)$.

We set

$$\varphi_{01m}(x_1) = \mu_m \cos(\alpha_m x_1), \quad m = 0,1,2,\cdots, \tag{25}$$

$$\varphi_{02m}(x_1) = \mu_m \sin(\alpha_m x_1), \quad m = 1,2,\cdots, \tag{26}$$

then for $k_1$ as a nonnegative even integer, we have

$$\varphi_{01m}^{(k_1)}(x_1) = \mu_m (-1)^{k_1/2} \alpha_m^{k_1} \cos(\alpha_m x_1), \tag{27}$$

$$\varphi_{02m}^{(k_1)}(x_1) = \mu_m(-1)^{k_1/2}\alpha_m^{k_1}\sin(\alpha_m x_1), \tag{28}$$

and for $k_1$ as a nonnegative odd integer, we have

$$\varphi_{01m}^{(k_1)}(x_1) = \mu_m(-1)^{(k_1+1)/2}\alpha_m^{k_1}\sin(\alpha_m x_1), \tag{29}$$

$$\varphi_{02m}^{(k_1)}(x_1) = \mu_m(-1)^{(k_1-1)/2}\alpha_m^{k_1}\cos(\alpha_m x_1). \tag{30}$$

We define the following vectors of basis functions as

$$\mathbf{\Phi}_{01}^{\mathrm{T}}(x_1) = [\varphi_{010}(x_1) \quad \varphi_{011}(x_1) \quad \cdots \quad \varphi_{01m}(x_1) \quad \cdots], \tag{31}$$

$$\mathbf{\Phi}_{02}^{\mathrm{T}}(x_1) = [\varphi_{021}(x_1) \quad \varphi_{022}(x_1) \quad \cdots \quad \varphi_{02m}(x_1) \quad \cdots], \tag{32}$$

$$\mathbf{\Phi}_{0}^{\mathrm{T}}(x_1) = [\mathbf{\Phi}_{01}^{\mathrm{T}}(x_1) \quad \mathbf{\Phi}_{02}^{\mathrm{T}}(x_1)], \tag{33}$$

and the corresponding vectors of derivatives of basis functions as

$$\mathbf{\Phi}_{01}^{(k_1)\mathrm{T}}(x_1) = [\varphi_{010}^{(k_1)}(x_1) \quad \varphi_{011}^{(k_1)}(x_1) \quad \cdots \quad \varphi_{01m}^{(k_1)}(x_1) \quad \cdots], \tag{34}$$

$$\mathbf{\Phi}_{02}^{(k_1)\mathrm{T}}(x_1) = [\varphi_{021}^{(k_1)}(x_1) \quad \varphi_{022}^{(k_1)}(x_1) \quad \cdots \quad \varphi_{02m}^{(k_1)}(x_1) \quad \cdots], \tag{35}$$

$$\mathbf{\Phi}_{0}^{(k_1)\mathrm{T}}(x_1) = [\mathbf{\Phi}_{01}^{(k_1)\mathrm{T}}(x_1) \quad \mathbf{\Phi}_{02}^{(k_1)\mathrm{T}}(x_1)]. \tag{36}$$

In addition, we define the following vectors of Fourier coefficients of $\varphi_0(x_1)$ as

$$\mathbf{q}_{01}^{\mathrm{T}} = [V_{010} \quad V_{011} \quad \cdots \quad V_{01m} \quad \cdots], \tag{37}$$

$$\mathbf{q}_{02}^{\mathrm{T}} = [V_{021} \quad V_{022} \quad \cdots \quad V_{02m} \quad \cdots], \tag{38}$$

$$\mathbf{q}_{0}^{\mathrm{T}} = [\mathbf{q}_{01}^{\mathrm{T}} \quad \mathbf{q}_{02}^{\mathrm{T}}]. \tag{39}$$

Then the internal function $\varphi_0(x_1)$ can be reduced to

$$\varphi_0(x_1) = \mathbf{\Phi}_0^{\mathrm{T}}(x_1) \cdot \mathbf{q}_0. \tag{40}$$

3. Expressions of the function $u(x_1)$ and its derivatives

Since the internal function $\varphi_0(x_1)$ satisfies the sufficient conditions for $2r$ times term-by-term differentiation of the full-range Fourier series of one-dimensional functions, then for any integer $k_1 = 0, 1, \cdots, 2r$, we have

$$\begin{aligned} u^{(k_1)}(x_1) &= \varphi_0^{(k_1)}(x_1) + \varphi_1^{(k_1)}(x_1) \\ &= \mathbf{\Phi}_0^{(k_1)\mathrm{T}}(x_1) \cdot \mathbf{q}_0 + \mathbf{\Phi}_1^{(k_1)\mathrm{T}}(x_1) \cdot \mathbf{q}_1. \end{aligned} \tag{41}$$

*3.2. For two-dimensional functions*

As to the first situation in Section 2.2 (similarly for other four situations), we expand the two-dimensional function $u(x_1, x_2)$ and its partial derivatives in full-range composite Fourier series.

1. Constructing the corner function $\varphi_3(x_1, x_2)$

We define the set

$$\bar{\Omega}_3 = \{(j,l) \mid j, l = 0, 1, 2, \cdots, j + l \leq 2r - 2\}, \tag{42}$$

and select the vector of supplementary functions as follows

$$\mathbf{p}_3^{\mathrm{T}}(x_1, x_2) = [p_{3,00}(x_1, x_2) \quad \cdots \quad p_{3,jl}(x_1, x_2) \quad \cdots \quad p_{3,0(2r-2)}(x_1, x_2)]_{(j,l)\in\bar{\Omega}_3}, \tag{43}$$

where $p_{3,jl}(x_1, x_2)$, $(j,l) \in \bar{\Omega}_3$, represent a set of $2r$ times continuously differentiable functions defined over the domain $[-a,a] \times [-b,b]$.

Then for $(k_1, k_2) \in \bar{\Omega}_3$, we define the corresponding vector of partial derivatives of the supplementary functions as
$$\mathbf{p}_3^{(k_1,k_2)\mathrm{T}}(x_1,x_2) = [p_{3,00}^{(k_1,k_2)}(x_1,x_2) \quad \cdots \quad p_{3,jl}^{(k_1,k_2)}(x_1,x_2) \quad \cdots \quad p_{3,0(2r-2)}^{(k_1,k_2)}(x_1,x_2)]_{(j,l)\in\bar{\Omega}_3}. \tag{44}$$

Further, let the vector of undetermined constants
$$\mathbf{a}_3^\mathrm{T} = [G_{3,00} \quad \cdots \quad G_{3,jl} \quad \cdots \quad G_{3,0(2r-2)}]_{(j,l)\in\bar{\Omega}_3}; \tag{45}$$

thus, we can write the corner function in the form
$$\varphi_3(x_1,x_2) = \mathbf{p}_3^\mathrm{T}(x_1,x_2)\cdot\mathbf{a}_3. \tag{46}$$

Combining the necessary conditions as given in Eq. (7) for the corner function $\varphi_3(x_1,x_2)$, we have
$$\mathbf{R}_3 \mathbf{a}_3 = \mathbf{q}_3, \tag{47}$$

where
$$\mathbf{R}_3 = \begin{bmatrix} \mathbf{p}_3^{(0,0)\mathrm{T}}(a,b) - \mathbf{p}_3^{(0,0)\mathrm{T}}(a,-b) - \mathbf{p}_3^{(0,0)\mathrm{T}}(-a,b) + \mathbf{p}_3^{(0,0)\mathrm{T}}(-a,-b) \\ \vdots \\ \mathbf{p}_3^{(k_1,k_2)\mathrm{T}}(a,b) - \mathbf{p}_3^{(k_1,k_2)\mathrm{T}}(a,-b) - \mathbf{p}_3^{(k_1,k_2)\mathrm{T}}(-a,b) + \mathbf{p}_3^{(k_1,k_2)\mathrm{T}}(-a,-b) \\ \vdots \\ \mathbf{p}_3^{(0,2r-2)\mathrm{T}}(a,b) - \mathbf{p}_3^{(0,2r-2)\mathrm{T}}(a,-b) - \mathbf{p}_3^{(0,2r-2)\mathrm{T}}(-a,b) + \mathbf{p}_3^{(0,2r-2)\mathrm{T}}(-a,-b) \end{bmatrix}_{(k_1,k_2)\in\bar{\Omega}_3} \tag{48}$$

and the vector involving corner values of the function $u(x_1,x_2)$ and its partial derivatives
$$\begin{aligned}\mathbf{q}_3^\mathrm{T} = [&u^{(0,0)}(a,b) - u^{(0,0)}(a,-b) - u^{(0,0)}(-a,b) + u^{(0,0)}(-a,-b) \quad \cdots \\ &u^{(k_1,k_2)}(a,b) - u^{(k_1,k_2)}(a,-b) - u^{(k_1,k_2)}(-a,b) + u^{(k_1,k_2)}(-a,-b) \quad \cdots \\ &u^{(0,2r-2)}(a,b) - u^{(0,2r-2)}(a,-b) - u^{(0,2r-2)}(-a,b) + u^{(0,2r-2)}(-a,-b)]_{(k_1,k_2)\in\bar{\Omega}_3}. \end{aligned} \tag{49}$$

Then we have
$$\mathbf{a}_3 = \mathbf{R}_3^{-1}\mathbf{q}_3. \tag{50}$$

If we define the vector of basis functions as
$$\mathbf{\Phi}_3^\mathrm{T}(x_1,x_2) = \mathbf{p}_3^\mathrm{T}(x_1,x_2)\cdot\mathbf{R}_3^{-1}, \tag{51}$$

then the corner function $\varphi_3(x_1,x_2)$ can be rewritten in matrix form as
$$\varphi_3(x_1,x_2) = \mathbf{\Phi}_3^\mathrm{T}(x_1,x_2)\cdot\mathbf{q}_3. \tag{52}$$

2. Constructing the boundary function $\varphi_1(x_1,x_2)$

Expand $\varphi_1(x_1,x_2)$ in a one-dimensional full-range Fourier series over the interval $[-b,b]$ along the $x_2$-direction
$$\varphi_1(x_1,x_2) = \sum_{n=0}^{\infty}\mu_n[\xi_{1n}(x_1)\cos(\beta_n x_2) + \xi_{2n}(x_1)\sin(\beta_n x_2)], \tag{53}$$

where $\beta_n = n\pi/b$, $\mu_n = \begin{cases} 1/2, & n=0 \\ 1, & n>0 \end{cases}$, $\xi_{1n}(x_1)$ and $\xi_{2n}(x_1)$ are the corresponding one-dimensional full-range Fourier coefficients.

To ensure the satisfaction of the necessary conditions for the boundary function $\varphi_1(x_1,x_2)$, we substitute the assumed series expansion, eq. (53), in Eqs. (5)-(7). We find that Eqs. (6)-(7) are automatically satisfied, and eq. (5) results in
$$\xi_{1n}^{(k_1)}(a) - \xi_{1n}^{(k_1)}(-a) = c_{1n}^{(k_1)} - d_{1n}^{(k_1)}, \quad n=0,1,2,\cdots, \quad k_1=0,1,\cdots,2r-1, \tag{54}$$
$$\xi_{2n}^{(k_1)}(a) - \xi_{2n}^{(k_1)}(-a) = c_{2n}^{(k_1)} - d_{2n}^{(k_1)}, \quad n=1,2,\cdots, \quad k_1=0,1,\cdots,2r-1, \tag{55}$$

where $c_{1n}^{(k_1)}$, $c_{2n}^{(k_1)}$, $d_{1n}^{(k_1)}$ and $d_{2n}^{(k_1)}$ are the Fourier coefficients of functions

$u^{(k_1,0)}(a,x_2) - \varphi_3^{(k_1,0)}(a,x_2)$ and $u^{(k_1,0)}(-a,x_2) - \varphi_3^{(k_1,0)}(-a,x_2)$ respectively, or we call them the boundary Fourier coefficients (at the opposite edges $x_1 = -a$ and $x_1 = a$) of the function $u(x_1,x_2) - \varphi_3(x_1,x_2)$ (see section 3.2 in [60]).

We select the vector of supplementary functions as follows
$$\mathbf{p}_{1n}^{\mathrm{T}}(x_1) = [p_{1n,1}(x_1) \quad p_{1n,2}(x_1) \quad \cdots \quad p_{1n,2r}(x_1)], \quad n=0,1,2,\cdots, \tag{56}$$
where $p_{1n,j}(x_1)$, $j=1,2,\cdots,2r$, are a set of $2r$ times continuously differentiable functions defined over the interval $[-a,a]$.

Then for $k_1 = 0,1,\cdots,2r$, we define the vector of derivatives of the selected functions as
$$\mathbf{p}_{1n}^{(k_1)\mathrm{T}}(x_1) = [p_{1n,1}^{(k_1)}(x_1) \quad p_{1n,2}^{(k_1)}(x_1) \quad \cdots \quad p_{1n,2r}^{(k_1)}(x_1)]. \tag{57}$$

Further, let the vectors of undetermined constants
$$\mathbf{a}_{1n,1}^{\mathrm{T}} = [G_{1n,1}^1 \quad G_{1n,1}^2 \quad \cdots \quad G_{1n,1}^{2r}], \quad n=0,1,2,\cdots, \tag{58}$$
$$\mathbf{a}_{1n,2}^{\mathrm{T}} = [G_{1n,2}^1 \quad G_{1n,2}^2 \quad \cdots \quad G_{1n,2}^{2r}], \quad n=1,2,\cdots, \tag{59}$$
thus we can construct the functions $\xi_{1n}(x_1)$ and $\xi_{2n}(x_1)$ as follows
$$\xi_{1n}(x_1) = \mathbf{p}_{1n}^{\mathrm{T}}(x_1) \cdot \mathbf{a}_{1n,1}, \quad n=0,1,2,\cdots, \tag{60}$$
$$\xi_{2n}(x_1) = \mathbf{p}_{1n}^{\mathrm{T}}(x_1) \cdot \mathbf{a}_{1n,2}, \quad n=1,2,\cdots. \tag{61}$$

Substitute the above two expressions into Eqs. (54) and (55) respectively, and we obtain
$$\mathbf{R}_{1n}\mathbf{a}_{1n,1} = \mathbf{q}_{1n,1}, \quad n=0,1,2,\cdots, \tag{62}$$
$$\mathbf{R}_{1n}\mathbf{a}_{1n,2} = \mathbf{q}_{1n,2}, \quad n=1,2,\cdots, \tag{63}$$
where
$$\mathbf{R}_{1n} = \begin{bmatrix} \mathbf{p}_{1n}^{(0)\mathrm{T}}(a) - \mathbf{p}_{1n}^{(0)\mathrm{T}}(-a) \\ \mathbf{p}_{1n}^{(1)\mathrm{T}}(a) - \mathbf{p}_{1n}^{(1)\mathrm{T}}(-a) \\ \vdots \\ \mathbf{p}_{1n}^{(2r-1)\mathrm{T}}(a) - \mathbf{p}_{1n}^{(2r-1)\mathrm{T}}(-a) \end{bmatrix}, \tag{64}$$
and the sub-vectors of boundary Fourier coefficients of the function $u(x_1,x_2) - \varphi_3(x_1,x_2)$
$$\mathbf{q}_{1n,1}^{\mathrm{T}} = [c_{1n}^{(0)} - d_{1n}^{(0)} \quad c_{1n}^{(1)} - d_{1n}^{(1)} \quad \cdots \quad c_{1n}^{(2r-1)} - d_{1n}^{(2r-1)}], \tag{65}$$
$$\mathbf{q}_{1n,2}^{\mathrm{T}} = [c_{2n}^{(0)} - d_{2n}^{(0)} \quad c_{2n}^{(1)} - d_{2n}^{(1)} \quad \cdots \quad c_{2n}^{(2r-1)} - d_{2n}^{(2r-1)}]. \tag{66}$$

Then we have
$$\mathbf{a}_{1n,1} = \mathbf{R}_{1n}^{-1}\mathbf{q}_{1n,1}, \quad n=0,1,2,\cdots, \tag{67}$$
$$\mathbf{a}_{1n,2} = \mathbf{R}_{1n}^{-1}\mathbf{q}_{1n,2}, \quad n=1,2,\cdots. \tag{68}$$

We define the sub-vectors of basis functions as
$$\mathbf{\Phi}_{1n,1}^{\mathrm{T}}(x_1,x_2) = \mu_n \cos(\beta_n x_2) \cdot \mathbf{p}_{1n}^{\mathrm{T}}(x_1) \cdot \mathbf{R}_{1n}^{-1}, \quad n=0,1,2,\cdots, \tag{69}$$
$$\mathbf{\Phi}_{1n,2}^{\mathrm{T}}(x_1,x_2) = \mu_n \sin(\beta_n x_2) \cdot \mathbf{p}_{1n}^{\mathrm{T}}(x_1) \cdot \mathbf{R}_{1n}^{-1}, \quad n=1,2,\cdots. \tag{70}$$

Then for $k_1 = 0,2,\cdots,2r$, we have
$$\mathbf{\Phi}_{1n,1}^{(k_1,k_2)\mathrm{T}}(x_1,x_2) = \mu_n (-1)^{k_2/2} \beta_n^{k_2} \cos(\beta_n x_2) \cdot \mathbf{p}_{1n}^{(k_1)\mathrm{T}}(x_1) \cdot \mathbf{R}_{1n}^{-1}, \tag{71}$$
$$\mathbf{\Phi}_{1n,2}^{(k_1,k_2)\mathrm{T}}(x_1,x_2) = \mu_n (-1)^{k_2/2} \beta_n^{k_2} \sin(\beta_n x_2) \cdot \mathbf{p}_{1n}^{(k_1)\mathrm{T}}(x_1) \cdot \mathbf{R}_{1n}^{-1}, \tag{72}$$
and for $k_1 = 1,3,\cdots,2r-1$, we have
$$\mathbf{\Phi}_{1n,1}^{(k_1,k_2)\mathrm{T}}(x_1,x_2) = \mu_n (-1)^{(k_2+1)/2} \beta_n^{k_2} \sin(\beta_n x_2) \cdot \mathbf{p}_{1n}^{(k_1)\mathrm{T}}(x_1) \cdot \mathbf{R}_{1n}^{-1}, \tag{73}$$
$$\mathbf{\Phi}_{1n,2}^{(k_1,k_2)\mathrm{T}}(x_1,x_2) = \mu_n (-1)^{(k_2-1)/2} \beta_n^{k_2} \cos(\beta_n x_2) \cdot \mathbf{p}_{1n}^{(k_1)\mathrm{T}}(x_1) \cdot \mathbf{R}_{1n}^{-1}. \tag{74}$$

We define the following vectors of basis functions as

$$\boldsymbol{\Phi}_{1,1}^{\mathrm{T}}(x_1,x_2) = [\boldsymbol{\Phi}_{10,1}^{\mathrm{T}}(x_1,x_2) \quad \boldsymbol{\Phi}_{11,1}^{\mathrm{T}}(x_1,x_2) \quad \cdots \quad \boldsymbol{\Phi}_{1n,1}^{\mathrm{T}}(x_1,x_2) \quad \cdots], \tag{75}$$

$$\boldsymbol{\Phi}_{1,2}^{\mathrm{T}}(x_1,x_2) = [\boldsymbol{\Phi}_{11,2}^{\mathrm{T}}(x_1,x_2) \quad \boldsymbol{\Phi}_{12,2}^{\mathrm{T}}(x_1,x_2) \quad \cdots \quad \boldsymbol{\Phi}_{1n,2}^{\mathrm{T}}(x_1,x_2) \quad \cdots], \tag{76}$$

$$\boldsymbol{\Phi}_{1}^{\mathrm{T}}(x_1,x_2) = [\boldsymbol{\Phi}_{1,1}^{\mathrm{T}}(x_1,x_2) \quad \boldsymbol{\Phi}_{1,2}^{\mathrm{T}}(x_1,x_2)], \tag{77}$$

and the corresponding vectors of partial derivatives of basis functions as

$$\boldsymbol{\Phi}_{1,1}^{(k_1,k_2)\mathrm{T}}(x_1,x_2) = [\boldsymbol{\Phi}_{10,1}^{(k_1,k_2)\mathrm{T}}(x_1,x_2) \quad \boldsymbol{\Phi}_{11,1}^{(k_1,k_2)\mathrm{T}}(x_1,x_2) \quad \cdots \quad \boldsymbol{\Phi}_{1n,1}^{(k_1,k_2)\mathrm{T}}(x_1,x_2) \quad \cdots], \tag{78}$$

$$\boldsymbol{\Phi}_{1,2}^{(k_1,k_2)\mathrm{T}}(x_1,x_2) = [\boldsymbol{\Phi}_{11,2}^{(k_1,k_2)\mathrm{T}}(x_1,x_2) \quad \boldsymbol{\Phi}_{12,2}^{(k_1,k_2)\mathrm{T}}(x_1,x_2) \quad \cdots \quad \boldsymbol{\Phi}_{1n,2}^{(k_1,k_2)\mathrm{T}}(x_1,x_2) \quad \cdots], \tag{79}$$

$$\boldsymbol{\Phi}_{1}^{(k_1,k_2)\mathrm{T}}(x_1,x_2) = [\boldsymbol{\Phi}_{1,1}^{(k_1,k_2)\mathrm{T}}(x_1,x_2) \quad \boldsymbol{\Phi}_{1,2}^{(k_1,k_2)\mathrm{T}}(x_1,x_2)]. \tag{80}$$

Meanwhile, we define the vectors of boundary Fourier coefficients of the function $u(x_1,x_2) - \varphi_3(x_1,x_2)$ as

$$\mathbf{q}_{1,1}^{\mathrm{T}} = [\mathbf{q}_{10,1}^{\mathrm{T}} \quad \mathbf{q}_{11,1}^{\mathrm{T}} \quad \cdots \quad \mathbf{q}_{1n,1}^{\mathrm{T}} \quad \cdots], \tag{81}$$

$$\mathbf{q}_{1,2}^{\mathrm{T}} = [\mathbf{q}_{11,2}^{\mathrm{T}} \quad \mathbf{q}_{12,2}^{\mathrm{T}} \quad \cdots \quad \mathbf{q}_{1n,2}^{\mathrm{T}} \quad \cdots], \tag{82}$$

$$\mathbf{q}_{1}^{\mathrm{T}} = [\mathbf{q}_{1,1}^{\mathrm{T}} \quad \mathbf{q}_{1,2}^{\mathrm{T}}], \tag{83}$$

then the boundary function $\varphi_1(x_1,x_2)$ can be obtained as

$$\varphi_1(x_1,x_2) = \boldsymbol{\Phi}_{1}^{\mathrm{T}}(x_1,x_2) \cdot \mathbf{q}_1. \tag{84}$$

3. Constructing the boundary function $\varphi_2(x_1,x_2)$

Similarly, we expand $\varphi_2(x_1,x_2)$ in a one-dimensional full-range Fourier series over the interval $[-a,a]$ along the $x_1$-direction

$$\varphi_2(x_1,x_2) = \sum_{m=0}^{\infty} \mu_m [\zeta_{1m}(x_2)\cos(\alpha_m x_1) + \zeta_{2m}(x_2)\sin(\alpha_m x_1)], \tag{85}$$

where $\alpha_m = m\pi/a$, $\mu_m = \begin{cases} 1/2, & m=0 \\ 1, & m>0 \end{cases}$, $\zeta_{1m}(x_2)$ and $\zeta_{2m}(x_2)$ are the corresponding one-dimensional full-range Fourier coefficients.

We denote by $p_{2m,l}(x_2)$, $l=1,2,\cdots,2r$, the selected supplementary functions that are $2r$ times continuously differentiable over the interval $[-b,b]$. By the same procedure, we obtain the expression of the vector of basis functions $\boldsymbol{\Phi}_2(x_1,x_2)$. And accordingly, the expression of the boundary function $\varphi_2(x_1,x_2)$ is derived as

$$\varphi_2(x_1,x_2) = \boldsymbol{\Phi}_{2}^{\mathrm{T}}(x_1,x_2) \cdot \mathbf{q}_2, \tag{86}$$

where $\mathbf{q}_2$ is the vector of boundary Fourier coefficients (at the opposite edges $x_2 = -b$ and $x_2 = b$) of the function $u(x_1,x_2) - \varphi_3(x_1,x_2)$.

4. Constructing the internal function $\varphi_0(x_1,x_2)$

Expand $\varphi_0(x_1,x_2)$ in a full-range Fourier series over the domain $[-a,a] \times [-b,b]$, then we have

$$\varphi_0(x_1,x_2) = \sum_{m=0}^{\infty}\sum_{n=0}^{\infty} \lambda_{mn} [V_{1mn}\cos(\alpha_m x_1)\cos(\beta_n x_2) + V_{2mn}\sin(\alpha_m x_1)\cos(\beta_n x_2)$$
$$+ V_{3mn}\cos(\alpha_m x_1)\sin(\beta_n x_2) + V_{4mn}\sin(\alpha_m x_1)\sin(\beta_n x_2)], \tag{87}$$

where $\alpha_m = m\pi/a$, $\beta_n = n\pi/b$,

$$\lambda_{mn} = \begin{cases} 1/4 & m=0, n=0 \\ 1/2 & m>0, n=0 \\ 1/2 & m=0, n>0 \\ 1 & m>0, n>0 \end{cases},$$

$V_{1mn}$, $V_{2mn}$, $V_{3mn}$ and $V_{4mn}$ are the Fourier coefficients of $\varphi_0(x_1, x_2)$ respectively.

We set

$$\varphi_{01mn}(x_1, x_2) = \lambda_{mn}\cos(\alpha_m x_1)\cos(\beta_n x_2), \quad m=0,1,2,\cdots, \quad n=0,1,2,\cdots, \tag{88}$$

$$\varphi_{02mn}(x_1, x_2) = \lambda_{mn}\sin(\alpha_m x_1)\cos(\beta_n x_2), \quad m=1,2,\cdots, \quad n=0,1,2,\cdots, \tag{89}$$

$$\varphi_{03mn}(x_1, x_2) = \lambda_{mn}\cos(\alpha_m x_1)\sin(\beta_n x_2), \quad m=0,1,2,\cdots, \quad n=1,2,\cdots, \tag{90}$$

$$\varphi_{04mn}(x_1, x_2) = \lambda_{mn}\sin(\alpha_m x_1)\sin(\beta_n x_2), \quad m=1,2,\cdots, \quad n=1,2,\cdots. \tag{91}$$

If the nonnegative integers $k_1$ and $k_2$ are even numbers, we have

$$\varphi_{01mn}^{(k_1,k_2)}(x_1, x_2) = \lambda_{mn}(-1)^{(k_1+k_2)/2}\alpha_m^{k_1}\beta_n^{k_2}\cos(\alpha_m x_1)\cos(\beta_n x_2), \tag{92}$$

$$\varphi_{02mn}^{(k_1,k_2)}(x_1, x_2) = \lambda_{mn}(-1)^{(k_1+k_2)/2}\alpha_m^{k_1}\beta_n^{k_2}\sin(\alpha_m x_1)\cos(\beta_n x_2), \tag{93}$$

$$\varphi_{03mn}^{(k_1,k_2)}(x_1, x_2) = \lambda_{mn}(-1)^{(k_1+k_2)/2}\alpha_m^{k_1}\beta_n^{k_2}\cos(\alpha_m x_1)\sin(\beta_n x_2), \tag{94}$$

$$\varphi_{04mn}^{(k_1,k_2)}(x_1, x_2) = \lambda_{mn}(-1)^{(k_1+k_2)/2}\alpha_m^{k_1}\beta_n^{k_2}\sin(\alpha_m x_1)\sin(\beta_n x_2). \tag{95}$$

And for other possible combinations of the nonnegative integers $k_1$ and $k_2$, similar equations are readily obtained.

We define the following vectors of basis functions as

$$\boldsymbol{\Phi}_{01}^{\mathrm{T}}(x_1, x_2) = [\varphi_{0100}(x_1, x_2) \quad \cdots \quad \varphi_{01mn}(x_1, x_2) \quad \cdots], \tag{96}$$

$$\boldsymbol{\Phi}_{02}^{\mathrm{T}}(x_1, x_2) = [\varphi_{0210}(x_1, x_2) \quad \cdots \quad \varphi_{02mn}(x_1, x_2) \quad \cdots], \tag{97}$$

$$\boldsymbol{\Phi}_{03}^{\mathrm{T}}(x_1, x_2) = [\varphi_{0301}(x_1, x_2) \quad \cdots \quad \varphi_{03mn}(x_1, x_2) \quad \cdots], \tag{98}$$

$$\boldsymbol{\Phi}_{04}^{\mathrm{T}}(x_1, x_2) = [\varphi_{0411}(x_1, x_2) \quad \cdots \quad \varphi_{04mn}(x_1, x_2) \quad \cdots], \tag{99}$$

$$\boldsymbol{\Phi}_0^{\mathrm{T}}(x_1, x_2) = [\boldsymbol{\Phi}_{01}^{\mathrm{T}}(x_1, x_2) \quad \boldsymbol{\Phi}_{02}^{\mathrm{T}}(x_1, x_2) \quad \boldsymbol{\Phi}_{03}^{\mathrm{T}}(x_1, x_2) \quad \boldsymbol{\Phi}_{04}^{\mathrm{T}}(x_1, x_2)], \tag{100}$$

and the corresponding vectors of partial derivatives of basis functions as

$$\boldsymbol{\Phi}_{01}^{(k_1,k_2)\mathrm{T}}(x_1, x_2) = [\varphi_{0100}^{(k_1,k_2)}(x_1, x_2) \quad \cdots \quad \varphi_{01mn}^{(k_1,k_2)}(x_1, x_2) \quad \cdots], \tag{101}$$

$$\boldsymbol{\Phi}_{02}^{(k_1,k_2)\mathrm{T}}(x_1, x_2) = [\varphi_{0210}^{(k_1,k_2)}(x_1, x_2) \quad \cdots \quad \varphi_{02mn}^{(k_1,k_2)}(x_1, x_2) \quad \cdots], \tag{102}$$

$$\boldsymbol{\Phi}_{03}^{(k_1,k_2)\mathrm{T}}(x_1, x_2) = [\varphi_{0301}^{(k_1,k_2)}(x_1, x_2) \quad \cdots \quad \varphi_{03mn}^{(k_1,k_2)}(x_1, x_2) \quad \cdots], \tag{103}$$

$$\boldsymbol{\Phi}_{04}^{(k_1,k_2)\mathrm{T}}(x_1, x_2) = [\varphi_{0411}^{(k_1,k_2)}(x_1, x_2) \quad \cdots \quad \varphi_{04mn}^{(k_1,k_2)}(x_1, x_2) \quad \cdots], \tag{104}$$

$$\boldsymbol{\Phi}_0^{(k_1,k_2)\mathrm{T}}(x_1, x_2) = [\boldsymbol{\Phi}_{01}^{(k_1,k_2)\mathrm{T}}(x_1, x_2) \quad \boldsymbol{\Phi}_{02}^{(k_1,k_2)\mathrm{T}}(x_1, x_2) \quad \boldsymbol{\Phi}_{03}^{(k_1,k_2)\mathrm{T}}(x_1, x_2) \quad \boldsymbol{\Phi}_{04}^{(k_1,k_2)\mathrm{T}}(x_1, x_2)]. \tag{105}$$

In addition, we define the vectors of Fourier coefficients of $\varphi_0(x_1, x_2)$ as

$$\mathbf{q}_{01}^{\mathrm{T}} = [V_{0100} \quad \cdots \quad V_{01mn} \quad \cdots], \tag{106}$$

$$\mathbf{q}_{02}^{\mathrm{T}} = [V_{0210} \quad \cdots \quad V_{02mn} \quad \cdots], \tag{107}$$

$$\mathbf{q}_{03}^{\mathrm{T}} = [V_{0301} \quad \cdots \quad V_{03mn} \quad \cdots], \tag{108}$$

$$\mathbf{q}_{04}^{\mathrm{T}} = [V_{0411} \quad \cdots \quad V_{04mn} \quad \cdots], \tag{109}$$

$$\mathbf{q}_0^{\mathrm{T}} = [\mathbf{q}_{01}^{\mathrm{T}} \quad \mathbf{q}_{02}^{\mathrm{T}} \quad \mathbf{q}_{03}^{\mathrm{T}} \quad \mathbf{q}_{04}^{\mathrm{T}}], \tag{110}$$

then the internal function $\varphi_0(x_1, x_2)$ can be expressed as

$$\varphi_0(x_1, x_2) = \boldsymbol{\Phi}_0^{\mathrm{T}}(x_1, x_2) \cdot \mathbf{q}_0. \tag{111}$$

5. Expressions of the function $u(x_1, x_2)$ and its partial derivatives

Since the internal function $\varphi_0(x_1, x_2)$ satisfies the sufficient conditions for $2r$ times term-by-term differentiation of the full-range Fourier series expansion of two-dimensional functions, and the boundary function $\varphi_1(x_1, x_2)$ satisfies the sufficient conditions for $2r$ times term-by-term differentiation of the one-dimensional full-range Fourier series expansion along the $x_2$-direction, and the boundary function $\varphi_2(x_1, x_2)$ satisfies the sufficient conditions for $2r$ times term-by-term differentiation of the one-dimensional full-range Fourier series expansion along the $x_1$-direction, then for $k_1, k_2 = 0, 1, 2, \cdots$ and $k_1 + k_2 \leq 2r$, we have

$$\begin{aligned} u^{(k_1,k_2)}(x_1,x_2) &= \varphi_0^{(k_1,k_2)}(x_1,x_2) + \varphi_1^{(k_1,k_2)}(x_1,x_2) + \varphi_2^{(k_1,k_2)}(x_1,x_2) + \varphi_3^{(k_1,k_2)}(x_1,x_2) \\ &= \mathbf{\Phi}_0^{(k_1,k_2)\mathrm{T}}(x_1,x_2) \cdot \mathbf{q}_0 + \mathbf{\Phi}_1^{(k_1,k_2)\mathrm{T}}(x_1,x_2) \cdot \mathbf{q}_1 \\ &\quad + \mathbf{\Phi}_2^{(k_1,k_2)\mathrm{T}}(x_1,x_2) \cdot \mathbf{q}_2 + \mathbf{\Phi}_3^{(k_1,k_2)\mathrm{T}}(x_1,x_2) \cdot \mathbf{q}_3. \end{aligned} \tag{112}$$

## 4. Algebraical polynomial based composite Fourier series method

In section 3, we have established a framework of simultaneous approximation of functions and their (partial) derivatives by composite Fourier series. In the framework, $p_{1,j}(x_1)$, and $p_{3,jl}(x_1, x_2)$, $p_{1n,j}(x_1)$, $p_{2m,l}(x_2)$ are supplementary functions and can actually be any sufficiently smooth functions. This implies that there is a large (theoretically, an infinite) number of possible choices for these supplementary functions. For each specific choice, we are able to obtain a corresponding composite Fourier series expansion. Therefore, what is described in section 3 actually provides a general way of improving the accuracy and convergence of the Fourier series direct-expansion method.

As is shown in Table 1, algebraical polynomials are specifically employed for supplementary functions in the framework of the composite Fourier series method. This leads to the algebraical polynomial based composite Fourier series method.

Compared with the Fourier series direct-expansion method, the algebraical polynomial based composite Fourier series method involves various kinds of basis functions, including the functions in the complete trigonometrical system, the algebraical polynomials of specific degrees and the products of algebraical polynomials and trigonometrical functions. Accordingly, this method is not only capable of convergence in mean square error as the Fourier series direct-expansion method, but also has the reproducing property of complete algebraical polynomials of specific degrees:

1. With Fourier series expansion of the one-dimensional function specified as the full-range composite Fourier series over the interval $[-a, a]$ or half-range composite cosine series over the interval $[0, a]$, the algebraical polynomial based composite Fourier series method can reproduce a one-dimensional complete algebraical polynomial of degree $2r$.

2. With Fourier series expansion of the one-dimensional function specified as the half-range composite sine series over the interval $[0, a]$, the algebraical polynomial based composite Fourier series method can reproduce a one-dimensional complete algebraical polynomial of degree $2r - 1$.

3. With Fourier series expansion of the two-dimensional function specified as the full-range composite Fourier series over the domain $[-a, a] \times [-b, b]$, the algebraical polynomial based composite Fourier series method can reproduce a two-dimensional complete algebraical polynomial of degree $2r$.

4. With Fourier series expansion of the two-dimensional function specified as the half-range composite sine-sine series over the domain $[0,a]\times[0,b]$, the algebraical polynomial based composite Fourier series method can reproduce a two-dimensional complete algebraical polynomial of degree $2r-1$.

For the above four cases, the first, second and fourth ones are easily verified. And a detailed proof of the third case is deferred to Appendix A.

Table 1: Supplementary functions in the composite Fourier series method.

| No. | Composite Fourier series | Supplementary functions |
|---|---|---|
| 1 | Full-range composite Fourier series over the interval $[-a,a]$ | $p_{1,j}(x_1)=(x_1/a)^j$, $j=1,2,\cdots,2r$, <br> $\mathbf{p}_1^T(x_1)=[p_{1,1}(x_1) \quad p_{1,2}(x_1) \quad \cdots \quad p_{1,2r}(x_1)]$ |
| 2 | Half-range composite cosine series over the interval $[0,a]$ | $p_{1,j}(x_1)=(x_1/a)^j$, $j=1,2,\cdots,2r$, <br> $\mathbf{p}_1^T(x_1)=[p_{1,1}(x_1) \quad p_{1,2}(x_1) \quad \cdots \quad p_{1,2r}(x_1)]$ |
| 3 | Half-range composite sine series over the interval $[0,a]$ | $p_{1,j}(x_1)=(x_1/a)^{j-1}$, $j=1,2,\cdots,2r$, <br> $\mathbf{p}_1^T(x_1)=[p_{1,1}(x_1) \quad p_{1,2}(x_1) \quad \cdots \quad p_{1,2r}(x_1)]$ |
| 4 | Full-range composite Fourier series over the domain $[-a,a]\times[-b,b]$ | a. selection of $\mathbf{p}_3^T(x_1,x_2)$ <br> we denote the set <br> $\omega_{31}=\omega_{32}=\omega_{33}=\{(j,l)\mid j,l=0,1,2,\cdots,\ j+l\leq r-2\}$, <br> $\omega_{34}=\{(j,l)\mid j,l=0,1,2,\cdots,\ j+l\leq r-1\}$, <br> $p_{31,jl}(x_1,x_2)=(x_1/a)^{2j+2}(x_2/b)^{2l+2}$, $(j,l)\in\omega_{31}$, <br> $p_{32,jl}(x_1,x_2)=(x_1/a)^{2j+1}(x_2/b)^{2l+2}$, $(j,l)\in\omega_{32}$, <br> $p_{33,jl}(x_1,x_2)=(x_1/a)^{2j+2}(x_2/b)^{2l+1}$, $(j,l)\in\omega_{33}$, <br> $p_{34,jl}(x_1,x_2)=(x_1/a)^{2j+1}(x_2/b)^{2l+1}$, $(j,l)\in\omega_{34}$, <br> $\mathbf{p}_{31}^T(x_1,x_2)=[p_{31,00}(x_1,x_2) \quad \cdots \quad p_{31,jl}(x_1,x_2) \quad \cdots]_{(j,l)\in\omega_{31}}$, <br> $\mathbf{p}_{32}^T(x_1,x_2)=[p_{32,00}(x_1,x_2) \quad \cdots \quad p_{32,jl}(x_1,x_2) \quad \cdots]_{(j,l)\in\omega_{32}}$, <br> $\mathbf{p}_{33}^T(x_1,x_2)=[p_{33,00}(x_1,x_2) \quad \cdots \quad p_{33,jl}(x_1,x_2) \quad \cdots]_{(j,l)\in\omega_{33}}$, <br> $\mathbf{p}_{34}^T(x_1,x_2)=[p_{34,00}(x_1,x_2) \quad \cdots \quad p_{34,jl}(x_1,x_2) \quad \cdots]_{(j,l)\in\omega_{34}}$, <br> $\mathbf{p}_3^T(x_1,x_2)=[\mathbf{p}_{31}^T(x_1,x_2) \quad \mathbf{p}_{32}^T(x_1,x_2) \quad \mathbf{p}_{33}^T(x_1,x_2) \quad \mathbf{p}_{34}^T(x_1,x_2)]$ |

b. selection of $\mathbf{p}_{1n}^{T}(x_1)$, $n = 0,1,2,\cdots$

  $p_{1n,j}(x_1) = (x_1/a)^j$, $j = 1,2,\cdots,2r$,

  $\mathbf{p}_{1n}^{T}(x_1) = [p_{1n,1}(x_1) \quad p_{1n,2}(x_1) \quad \cdots \quad p_{1n,2r}(x_1)]$

  c. selection of $\mathbf{p}_{2m}^{T}(x_2)$, $m = 0,1,2,\cdots$

  $p_{2m,l}(x_2) = (x_2/b)^l$, $l = 1,2,\cdots,2r$,

  $\mathbf{p}_{2m}^{T}(x_2) = [p_{2m,1}(x_2) \quad p_{2m,2}(x_2) \quad \cdots \quad p_{2m,2r}(x_2)]$

| 5 | Half-range composite sine-sine series over the domain $[0,a] \times [0,b]$ | a. selection of $\mathbf{p}_3^T(x_1, x_2)$ ... |

5 | Half-range composite sine-sine series over the domain $[0,a]\times[0,b]$

a. selection of $\mathbf{p}_3^{T}(x_1, x_2)$

If $r$ is an even integer, we denote the set
$\Omega_3 = \{(j,l) | j, l = 0,1,2,\cdots, j+l \leq 2r-1;$
$\quad j + l = 2r, j = 1,2,\cdots,r/2 \text{ or } l = 1,2,\cdots,r/2\}$,

and if $r$ is an odd integer, we denote the set
$\Omega_3 = \{(j,l) | j, l = 0,1,2,\cdots, j+l \leq 2r-1; j+l = 2r,$
$j = 1,2,\cdots,(r-1)/2 \text{ or } l = 1,2,\cdots,(r-1)/2; j = l = r\}$,

$p_{3,jl}(x_1, x_2) = (x_1/a)^j (x_2/b)^l$, $(j,l) \in \Omega_3$,

$\mathbf{p}_3^{T}(x_1, x_2) = [p_{3,00}(x_1, x_2) \quad \cdots \quad p_{3,jl}(x_1, x_2) \quad \cdots \quad ]_{(j,l) \in \Omega_3}$

b. selection of $\mathbf{p}_{1n}^{T}(x_1)$, $n = 1,2,\cdots$

$p_{1n,j}(x_1) = (x_1/a)^{j-1}$, $j = 1,2,\cdots,2r$,

$\mathbf{p}_{1n}^{T}(x_1) = [p_{1n,1}(x_1) \quad p_{1n,2}(x_1) \quad \cdots \quad p_{1n,2r}(x_1)]$

c. selection of $\mathbf{p}_{2m}^{T}(x_2)$, $m = 1,2,\cdots$

$p_{2m,l}(x_2) = (x_2/b)^{l-1}$, $l = 1,2,\cdots,2r$,

$\mathbf{p}_{2m}^{T}(x_2) = [p_{2m,1}(x_2) \quad p_{2m,2}(x_2) \quad \cdots \quad p_{2m,2r}(x_2)]$

## 5. Numerical examples

In this section, with the computational parameter $2r = 6$ and length width ratio $a/b = 1.0$, we consider the four different kinds of sample functions as listed in Table 2, namely the one-dimensional algebraical polynomials, one-dimensional trigonometrical functions, two-dimensional algebraical polynomials and two-dimensional trigonometrical functions, and investigate the convergence characteristics and approximation accuracy of the algebraical polynomial based composite Fourier series method.

Table 2: Four different kinds of sample functions.

| No. | sample function |
|---|---|
| 1 | $u(x_1) = \dfrac{1}{2} - \dfrac{x_1}{a} - \dfrac{1}{2}(\dfrac{x_1}{a})^2 + (\dfrac{x_1}{a})^3$, $x_1 \in [-a, a]$ |
| 2 | $u(x_1) = \sin(\alpha_0 x_1)$, $\alpha_0 = \pi/2a$, $x_1 \in [-a, a]$ |
| 3 | $u(x_1) = \dfrac{1}{2} - \dfrac{x_1}{a} - \dfrac{1}{2}(\dfrac{x_1}{a})^2 + (\dfrac{x_1}{a})^3$, $x_1 \in [0, a]$ |

| | | |
|---|---|---|
| 4 | $u(x_1) = \cos(\alpha_0 x_1)$, $\alpha_0 = \pi/2a$, $x_1 \in [0,a]$ | |
| 5 | $u(x_1,x_2) = [\frac{1}{2} - \frac{x_1}{a} - \frac{1}{2}(\frac{x_1}{a})^2 + (\frac{x_1}{a})^3][\frac{1}{2} - \frac{x_2}{b} - \frac{1}{2}(\frac{x_2}{b})^2 + (\frac{x_2}{b})^3]$, $(x_1,x_2) \in [-a,a] \times [-b,b]$ | |
| 6 | $u(x_1,x_2) = \sin(\alpha_0 x_1)\cos(\beta_0 x_2)$, $\alpha_0 = \pi/2a$, $\beta_0 = \pi/2b$, $(x_1,x_2) \in [-a,a] \times [-b,b]$ | |
| 7 | $u(x_1,x_2) = [\frac{1}{2} - \frac{x_1}{a} - \frac{1}{2}(\frac{x_1}{a})^2 + (\frac{x_1}{a})^3][\frac{1}{2} - \frac{x_2}{b} - \frac{1}{2}(\frac{x_2}{b})^2 + (\frac{x_2}{b})^3]$, $(x_1,x_2) \in [0,a] \times [0,b]$ | |
| 8 | $u(x_1,x_2) = \sin(\alpha_0 x_1)\cos(\beta_0 x_2)$, $\alpha_0 = \pi/2a$, $\beta_0 = \pi/2b$, $(x_1,x_2) \in [0,a] \times [0,b]$ | |

*5.1. Error indexes of simultaneous approximation*

For the convenience of the application of the composite Fourier series method, we define some errors indexes of simultaneous approximation for describing difference of convergence speeds between different components of the (partial) derivative vectors and in different parts of the expansion domains.

Take the two-dimensional function $u(x_1, x_2)$ as an example. As shown in Figure 1.b, let $N_1^t$ and $N_2^t$ be positive integers and let $\Upsilon^t = \{(x_{1,n_1}, x_{2,n_2}), n_1 = 1, 2, \cdots, N_1^t, n_2 = 1, 2, \cdots, N_2^t\}$ be the set of uniformly distributed sampling points over the two-dimensional expansion domain. We conduct the decomposition of the set $\Upsilon^t$ into three subsets $\Upsilon_I^t$, $\Upsilon_B^t$ and $\Upsilon_C^t$, where $\Upsilon_I^t$ is the set of sampling points within the domain, $\Upsilon_B^t$ is the set of sampling points on the boundary (excluding the four corner points) and $\Upsilon_C^t$ is the set of sampling points at the corners. Then with the numbers of elements of $\Upsilon^t$ and its three subsets denoted by $I_\gamma^t$, $I_{\gamma I}^t$, $I_{\gamma B}^t$ and $I_{\gamma C}^t$ respectively, we have $I_\gamma^t = N_1^t N_2^t$, $I_{\gamma I}^t = (N_1^t - 2)(N_2^t - 2)$, $I_{\gamma B}^t = 2(N_1^t - 2) + 2(N_2^t - 2)$ and $I_{\gamma C}^t = 4$.

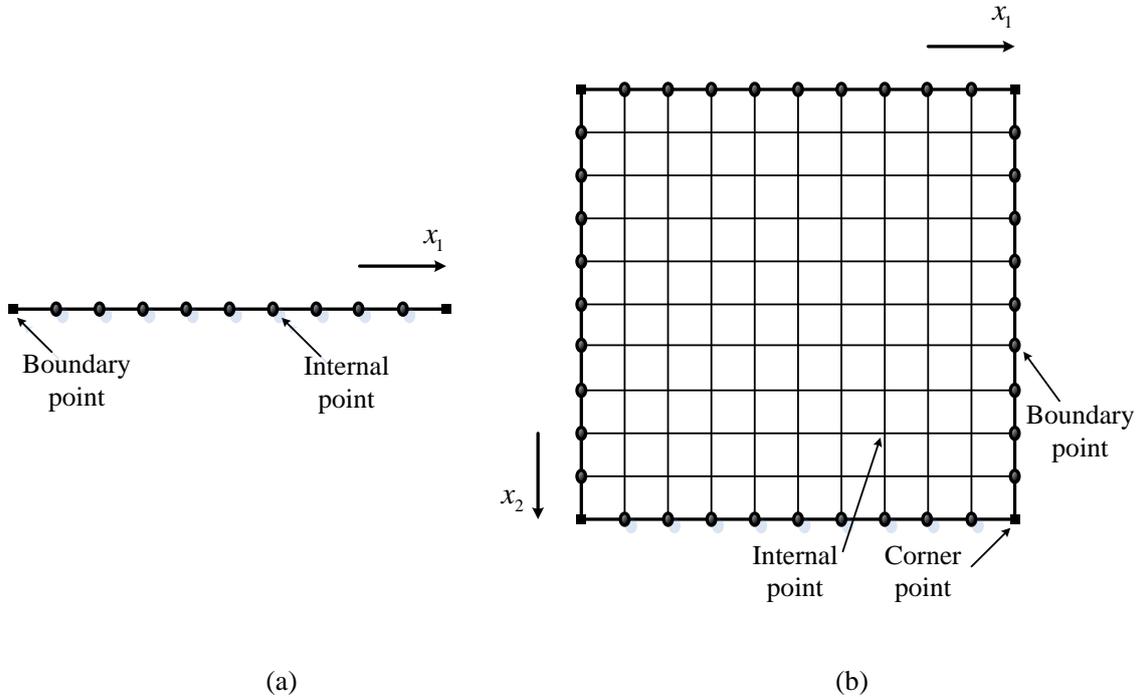

Figure 1: Sampling points in computation of errors in simultaneous composite Fourier series approximation: (a) for one-dimensional functions, (b) for two-dimensional functions.

For the two-dimensional function $u(x_1, x_2)$, let $M$ and $N$ be the numbers of truncated terms of the composite Fourier series respectively along the $x_1$- and $x_2$- direction, and we denote the corresponding partial sum by $u_{M,N}(x_1, x_2)$. Suppose that $k_1$, $k_2$ and $p$ are nonnegative integers, $k_1 + k_2 \leq 2r$ and $p \leq 2r$, then as presented in Table 3, we introduce some error indexes of simultaneous composite Fourier series approximation of the two-dimensional function $u(x_1, x_2)$ and its partial derivatives.

Table 3: Error indexes of simultaneous approximation of composite Fourier series for partial derivative vector of the two-dimensional function $u(x_1, x_2)$.

| | Set of sampling points | | | |
| --- | --- | --- | --- | --- |
| | $\Upsilon^t$ | $\Upsilon_I^t$ | $\Upsilon_B^t$ | $\Upsilon_C^t$ |
| Single component, $u^{(k_1,k_2)}(x_1,x_2)$, of partial derivative vector | $e^{(k_1,k_2)}(u_{M,N})$ | $e_I^{(k_1,k_2)}(u_{M,N})$ | $e_B^{(k_1,k_2)}(u_{M,N})$ | $e_C^{(k_1,k_2)}(u_{M,N})$ |
| $p$ th order components of partial derivative vector | $|e|^p(u_{M,N})$ | $|e_I|^p(u_{M,N})$ | $|e_B|^p(u_{M,N})$ | $|e_C|^p(u_{M,N})$ |
| Components up to $p$ th order of partial derivative vector | $\|e\|^p(u_{M,N})$ | $\|e_I\|^p(u_{M,N})$ | $\|e_B\|^p(u_{M,N})$ | $\|e_C\|^p(u_{M,N})$ |

In the error indexes of simultaneous approximation, the overall error of the single component, $u^{(k_1,k_2)}(x_1, x_2)$, of the partial derivative vector of the two-dimensional function $u(x_1, x_2)$ is defined by

$$e^{(k_1,k_2)}(u_{M,N}) = \frac{1}{I_\gamma^t \cdot u_{\max}^{(k_1,k_2)}} \sum_{(x_1,x_2) \in \Upsilon^t} \left| u_{M,N}^{(k_1,k_2)}(x_1,x_2) - u^{(k_1,k_2)}(x_1,x_2) \right|, \quad (113)$$

where $u_{\max}^{(k_1,k_2)}$ is the maximum value of the partial derivative $u^{(k_1,k_2)}(x_1, x_2)$ over the expansion domain. We also use $e(u_{M,N}^{(k_1,k_2)})$ or $e(u^{(k_1,k_2)})$ to denote $e^{(k_1,k_2)}(u_{M,N})$. Then on the special occasion $k_1 = k_2 = 0$, $e^{(k_1,k_2)}(u_{M,N})$ is further simplified as $e(u_{M,N})$ or $e(u)$.

The overall error of the $p$ th order components of the partial derivative vector of the two-dimensional function $u(x_1, x_2)$ is defined by

$$\begin{aligned} |e|^p(u_{M,N}) &= \frac{1}{p+1} \sum_{\substack{k_1 \geq 0, k_2 \geq 0 \\ k_1+k_2=p}} e^{(k_1,k_2)}(u_{M,N}) \\ &= \frac{1}{p+1} \sum_{\substack{k_1 \geq 0, k_2 \geq 0 \\ k_1+k_2=p}} \frac{1}{I_\gamma^t \cdot u_{\max}^{(k_1,k_2)}} \sum_{(x_1,x_2) \in \Upsilon^t} \left| u_{M,N}^{(k_1,k_2)}(x_1,x_2) - u^{(k_1,k_2)}(x_1,x_2) \right|. \end{aligned} \quad (114)$$

And the overall error of the components up to $p$ th order of the partial derivative vector of the two-dimensional function $u(x_1, x_2)$ is defined by

$$\|e\|^p(u_{M,N}) = \frac{2}{(p+1)(p+2)} \sum_{\substack{k_1 \geq 0, k_2 \geq 0 \\ k_1+k_2 \leq p}} e^{(k_1,k_2)}(u_{M,N})$$

$$= \frac{2}{(p+1)(p+2)} \sum_{\substack{k_1 \geq 0, k_2 \geq 0 \\ k_1+k_2 \leq p}} \frac{1}{I_\gamma^t \cdot u_{\max}^{(k_1,k_2)}} \sum_{(x_1,x_2) \in \Upsilon^t} \left| u_{M,N}^{(k_1,k_2)}(x_1,x_2) - u^{(k_1,k_2)}(x_1,x_2) \right|. \tag{115}$$

In Eqs. (113)-(115), by changing the set of sampling points $\Upsilon^t$ to $\Upsilon_I^t$, $\Upsilon_B^t$ and $\Upsilon_C^t$ respectively, we can obtain the internal error, boundary error and corner error of the single component, the $p$ th order components and the components up to $p$ th order of the partial derivative vector of the two-dimensional function $u(x_1,x_2)$.

With reference to Figure 1.a, we also introduce some error indexes of simultaneous composite Fourier series approximation of the one-dimensional function $u(x_1)$ and its derivatives (see Table 4), and the corresponding expressions are not presented here in the interest of brevity.

Table 4: Error indexes of simultaneous approximation of composite Fourier series for derivative vector of the one-dimensional function $u(x_1)$.

| | Set of sampling points | | |
|---|---|---|---|
| | $\Upsilon^t$ | $\Upsilon_I^t$ | $\Upsilon_B^t$ |
| Single component, $u^{(k_1)}(x_1)$, of derivative vector | $e^{(k_1)}(u_M)$ | $e_I^{(k_1)}(u_M)$ | $e_B^{(k_1)}(u_M)$ |
| Components up to $p$ th order of derivative vector | $\|e\|^p(u_M)$ | $\|e_I\|^p(u_M)$ | $\|e_B\|^p(u_M)$ |

*5.2. Convergence characteristics*

As to the sample functions given in Table 2, we truncate the composite Fourier series successively with the first 2, 3, 5, 10, 20, 30 and 40 terms. And then we set $N_1^t = N_2^t = 101$ and compute for $p = 0,1,\cdots,6$ the overall errors, internal errors, boundary errors, and corner errors of the single component, the $p$ th order component(s) and the components up to $p$ th order of the (partial) derivative vectors. Some of the results are shown in Figures 2-5.

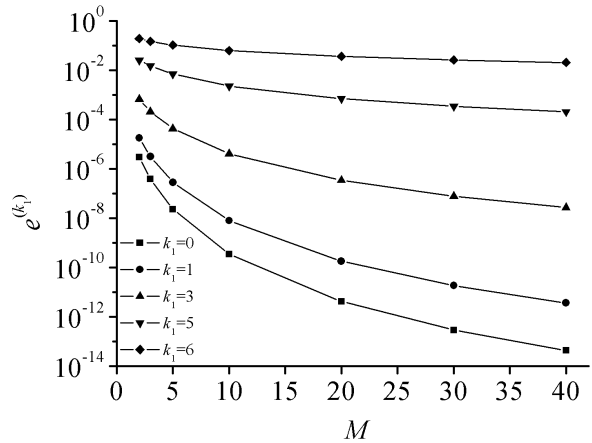

(a)

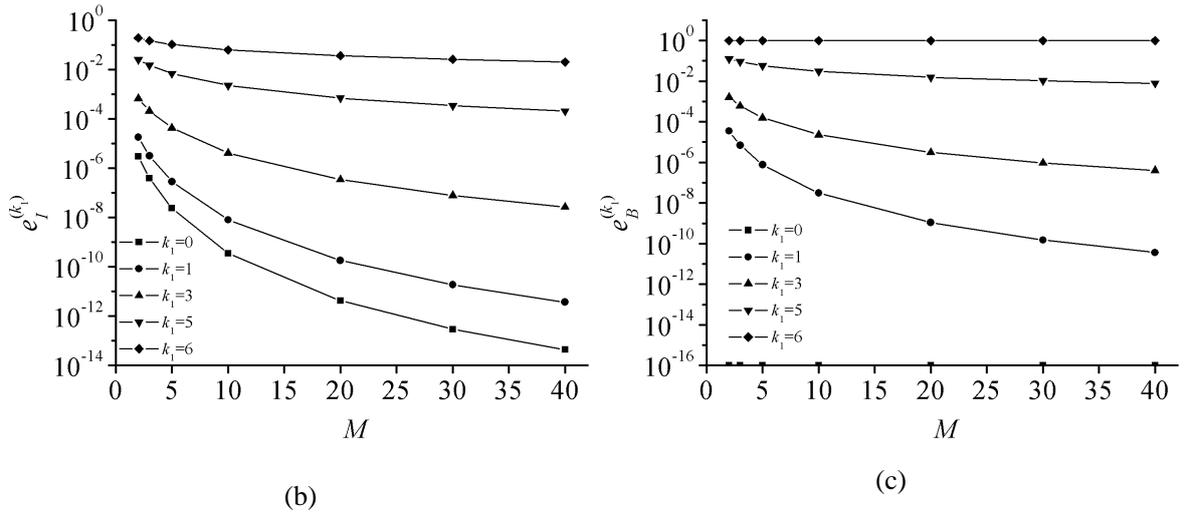

(b)              (c)

Figure 2: Convergence of the composite Fourier series of sample function 2 and its derivatives: (a) $e^{(k_1)}(u_M)$ -$M$ curves, (b) $e_I^{(k_1)}(u_M)$ -$M$ curves, (c) $e_B^{(k_1)}(u_M)$ -$M$ curves.

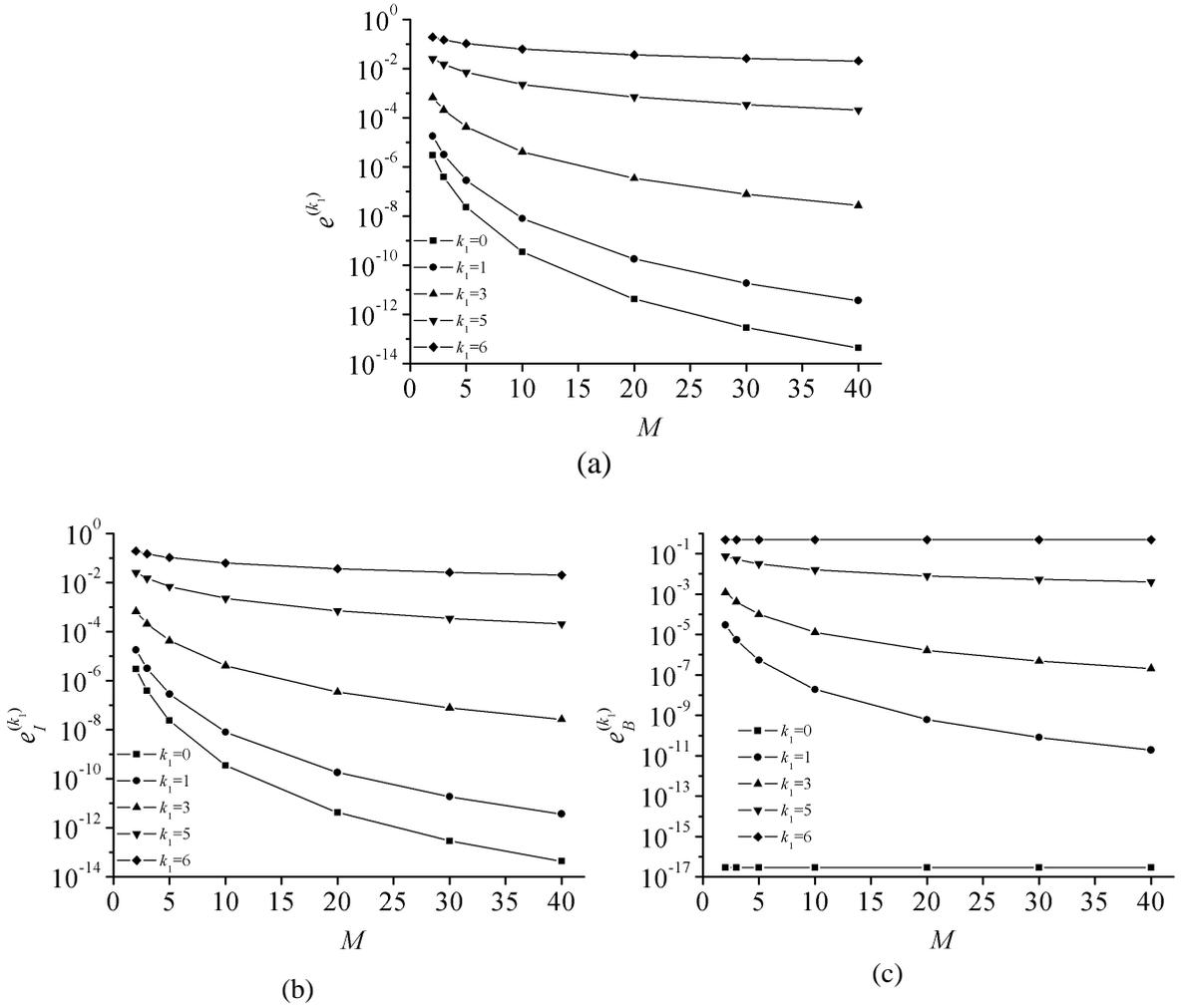

Figure 3: Convergence of the composite Fourier series of sample function 4 and its derivatives: (a) $e^{(k_1)}(u_M)$ -$M$ curves, (b) $e_I^{(k_1)}(u_M)$ -$M$ curves, (c) $e_B^{(k_1)}(u_M)$ -$M$ curves.

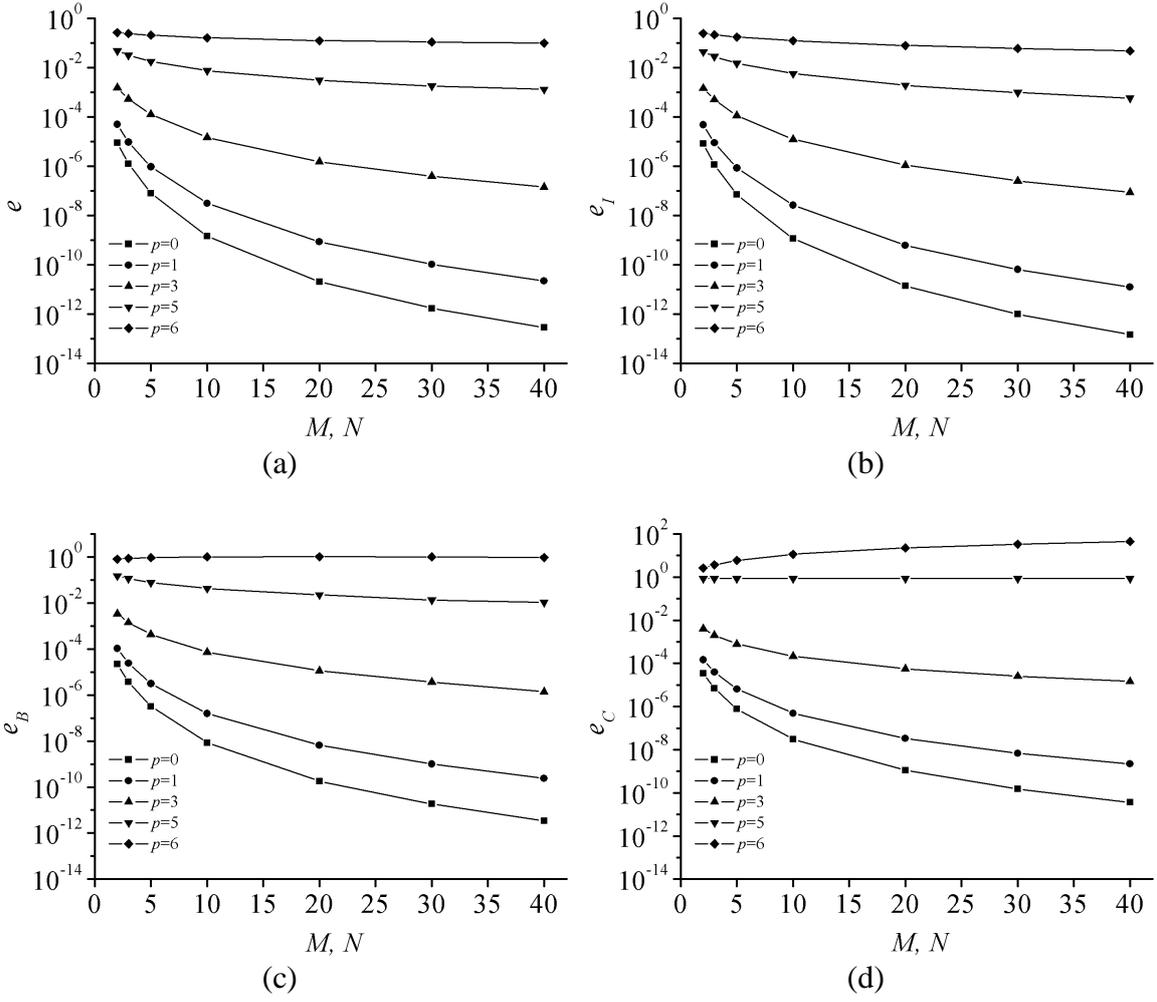

Figure 4: Convergence of the composite Fourier series of sample function 6 and its partial derivatives: (a) $|e|^p(u_{M,N})$ -$M, N$ curves, (b) $|e_I|^p(u_{M,N})$ -$M, N$ curves, (c) $|e_B|^p(u_{M,N})$ -$M, N$ curves, (d) $|e_C|^p(u_{M,N})$ -$M, N$ curves.

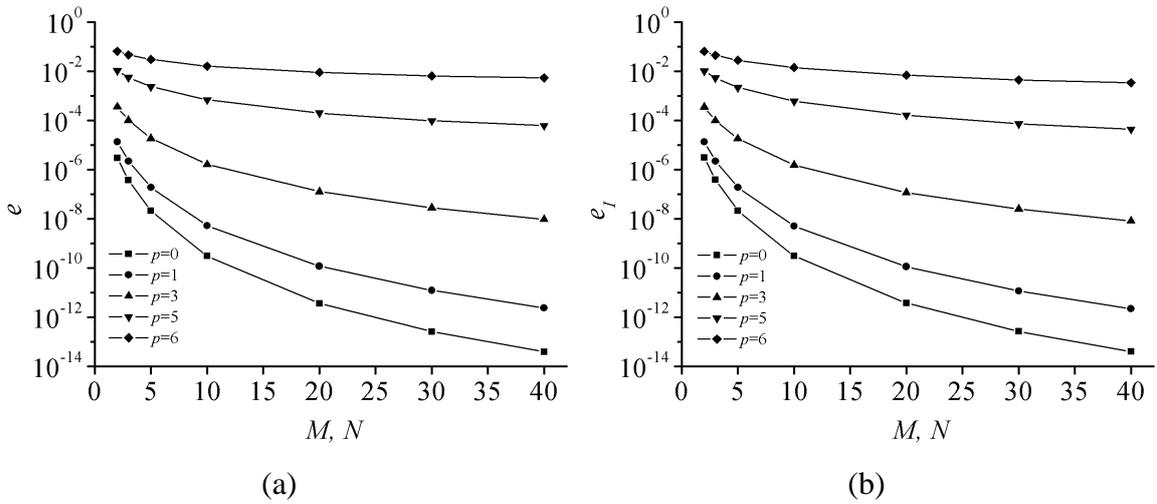

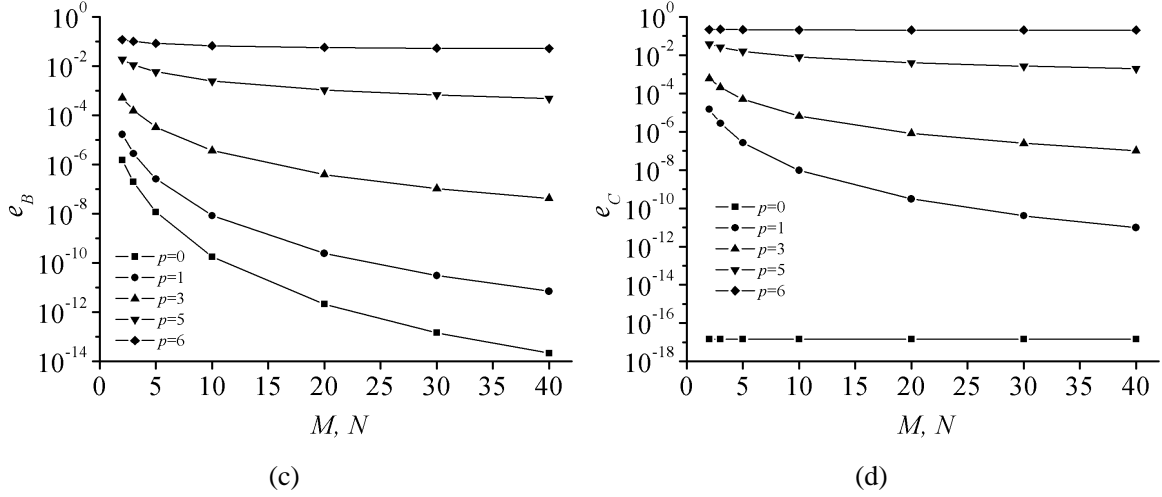

Figure 5: Convergence of the composite Fourier series of sample function 8 and its partial derivatives: (a) $|e|^p(u_{M,N})$ -$M, N$ curves, (b) $|e_I|^p(u_{M,N})$ -$M, N$ curves, (c) $|e_B|^p(u_{M,N})$ -$M, N$ curves, (d) $|e_C|^p(u_{M,N})$ -$M, N$ curves.

With the increase of the number of truncated terms of the composite Fourier series, the evolution of computed values of the errors is displayed as follows:

1. The overall errors, internal errors, boundary errors, and corner errors of the sample functions and their first to fourth order (partial) derivatives decrease rapidly as the number of truncated terms of the composite Fourier series increases. When the number of truncated terms equals to 2, the errors are respectively below 8.4E-3, 7.8E-3, 2.0E-2, 1.3E-1. When the number of truncated terms equals to 10, the errors are respectively below 3.3E-4, 2.7E-4, 1.6E-3, 3.0E-2, and show a sustained and rapid decrease. And when the number of truncated terms equals to 40, the errors are respectively below 1.5E-5, 7.2E-6, 1.4E-4, 7.9E-3.

2. The overall errors, internal errors, and boundary errors of the fifth order (partial) derivatives of the sample functions decrease rapidly as the number of truncated terms of the composite Fourier series increases. When the number of truncated terms equals to 2, the errors are respectively below 4.8E-2, 4.4E-2, 1.5E-1. When the number of truncated terms equals to 10, the errors are respectively below 7.5E-3, 5.7E-3, 4.3E-2, and show a sustained decrease. And when the number of truncated terms equals to 40, the errors are respectively below 1.3E-3, 5.8E-4, 1.0E-2.

3. The corner errors of the fifth order partial derivatives of the two-dimensional sample functions remain almost unchanged as the number of truncated terms of the composite Fourier series increases (see Figure 4.d) or show a slow decrease (see Figure 5.d). The errors are in the range of 2.0E-3 to 8.6E-1.

4. The overall errors and internal errors of the sixth order (partial) derivatives of the sample functions decrease slowly as the number of truncated terms of the composite Fourier series increases. When the number of truncated terms equals to 2, the errors are respectively below 2.6E-1, 2.4E-1. When the number of truncated terms equals to 10, the errors are respectively below 1.6E-1, 1.2E-1, and show a slow-down decrease. And when the number of truncated terms equals to 40, the errors are respectively below 1.0E-1, 4.7E-2.

5. The boundary errors of the sixth order (partial) derivatives of the sample function remain almost unchanged or show a slow decrease (see Figure 5.c) as the number of truncated terms of the composite Fourier series increases. The errors are between 5.2E-2 and 1.0E0.

6. The corner errors of the sixth order partial derivatives of the two-dimensional sample functions remain almost unchanged (see Figure 5.d) or show a sustained increase (see Figure 4.d) as the number of truncated terms of the composite Fourier series increases.

Therefore, the convergence characteristics of the algebraical polynomial based composite Fourier series method are summarized: the composite Fourier series of the function and its lower and middle order (partial) derivatives converge well; the composite Fourier series of the higher order (partial) derivatives of the function do not converge so well; the composite Fourier series of the function and its (partial) derivatives converge well within the expansion domain; the composite Fourier series of the function and its lower and middle order (partial) derivatives converge well on the boundary and at the corners of the expansion domain; the composite Fourier series of the higher order (partial) derivatives of the function do not converge so well on the boundary or at the corners of the expansion domain.

Meanwhile, the evolution of errors for the composite Fourier series approximation with the increasing number of truncated terms also sheds insight into the accuracy improvement of the algebraical polynomial based composite Fourier series. More precisely, as the number of truncated terms increases, the approximation accuracy of the algebraical polynomial based composite Fourier series is improved along the two dimensionalities, namely, the order of the (partial) derivatives of functions (from lower orders up to higher orders) and the expansion domain (gradual transition from the boundaries and corners to the internal part). To better understand this, we consider, for instance, sample function 8. Assume the overall error bound is specified, say, 1.0E-3. We find that the composite Fourier series with the first 2 terms has already yielded the results to this desired degree of accuracy for the sample function and its partial derivatives up to fourth order. In contrast, to obtain the results with the specified precision for the fifth order partial derivatives, the composite Fourier series is truncated to include the first 10 terms. And when the number of truncated terms increases to 20, and then 30, and finally to 40, the overall error for the sixth order partial derivatives gradually decreases to 5.4E-3 (still greater than the error bound). This implies that some more terms are needed in the calculation. On the other hand, when the number of truncated terms equals to 2, the boundary errors and corner errors of the sample function and its second and fourth order partial derivatives are less than their corresponding internal errors. This may be attributed to the introduction of the corner function and the boundary functions in the composite Fourier series. Since the discontinuities potentially to the sample function and its second and fourth order partial derivatives have been explicitly absorbed by the corner function and the boundary functions, the internal function in the form of the half-range sine-sine series now only represents a residual function that has at least fifth order continuous partial derivatives over the expansion domain. It is reflected in Figure 5 that the internal errors of the sample function and its partial derivatives up to sixth order decrease with higher speeds than their corresponding boundary errors and corner errors.

*5.3. Approximation accuracy*

In view of the excellent numerical behaviour of the algebraical polynomial based composite Fourier series, the settings $M = 40$ for one-dimensional sample functions and $M = N = 40$ for two-dimensional sample functions will be used in all the subsequent calculations.

1. Reproducing property of algebraical polynomials

Considering the sample functions in the form of algebraical polynomials, we present in Table 5 the computed results of the vectors of undetermined constants involved in the algebraical polynomial based composite Fourier series method. For the basis functions that are not algebraical polynomials, the orders of magnitude of the corresponding undetermined constants are all below 1.0E-13. And for the basis functions in the form of algebraical polynomials, the corresponding undetermined constants can be combined to precisely reconstruct the original sample functions. This verifies the reproducing property of algebraical polynomials for the algebraical polynomial based composite Fourier series method through numerical investigations.

It's worthy of notice that for the specified computational parameter $2r = 6$, the basis functions in the half-range composite sine-sine series expansion (over the domain $[0,a] \times [0,b]$) consist of not only the complete algebraical polynomials of degree 5, but also a algebraical polynomial, $(x_1/a)^3 (x_2/b)^3$, of degree 6. This explains the fact that the algebraical polynomial based composite Fourier series method can reproduce sample function 7, an incomplete two-dimensional algebraical polynomial of degree 6.

Table 5: Vectors of undetermined constants involved in sample functions in the form of algebraical polynomials.

| Sample functions | Vectors of undetermined constants |
| --- | --- |
| 1 | a. $\mathbf{q}_0^T$ <br> $\mathbf{q}_0^T = [\mathbf{q}_{01}^T \quad \mathbf{q}_{02}^T]$, <br> in which <br> $\mathbf{q}_{01}^T = [1.0 \quad 1.001396\text{E-}16 \quad 6.468621\text{E-}17 \quad \cdots \quad 7.378673\text{E-}17 \quad 1.233958\text{E-}17]$, <br> $\mathbf{q}_{02}^T = [8.012204\text{E-}18 \quad 1.184774\text{E-}16 \quad \cdots \quad -3.757968\text{E-}18 \quad 1.034444\text{E-}16]$ <br> b. $\mathbf{q}_1^T$ <br> $\mathbf{q}_1^T = [0.0 \quad -2.0 \quad 12.0 \quad 0.0 \quad 0.0 \quad 0.0]$ |
| 3 | a. $\mathbf{q}_{02}^T$ <br> $\mathbf{q}_{02}^T = [3.987988\text{E-}17 \quad 1.296269\text{E-}17 \quad \cdots \quad 2.783921\text{E-}19 \quad 4.874507\text{E-}19]$ <br> b. $\mathbf{q}_1^T$ <br> $\mathbf{q}_1^T = [0.0 \quad 5.0 \quad 0.0 \quad 0.5 \quad -1.0 \quad 0.0]$ |
| 5 | a. $\mathbf{q}_0^T$ <br> $\mathbf{q}_0^T = [\mathbf{q}_{01}^T \quad \mathbf{q}_{02}^T \quad \mathbf{q}_{03}^T \quad \mathbf{q}_{04}^T]$, <br> in which <br> $\mathbf{q}_{01}^T = [1.0 \quad 2.070504\text{E-}14 \quad -5.212722\text{E-}15 \quad \cdots \quad -2.718446\text{E-}19 \quad -2.488159\text{E-}20]$, <br> $\mathbf{q}_{02}^T = [3.052553\text{E-}13 \quad -1.604313\text{E-}17 \quad \cdots \quad -8.664684\text{E-}19 \quad 7.630576\text{E-}20]$, <br> $\mathbf{q}_{03}^T = [1.528638\text{E-}13 \quad -1.907658\text{E-}14 \quad \cdots \quad -3.189244\text{E-}19 \quad 3.043197\text{E-}19]$, <br> $\mathbf{q}_{04}^T = [-3.195236\text{E-}17 \quad -5.780299\text{E-}17 \quad \cdots \quad 1.158276\text{E-}19 \quad 3.719174\text{E-}19]$ <br> b. $\mathbf{q}_1^T$ <br> $\mathbf{q}_1^T = [\mathbf{q}_{11}^T \quad \mathbf{q}_{12}^T]$, <br> $\mathbf{q}_{11}^T = [\mathbf{q}_{10,1}^T \quad \mathbf{q}_{11,1}^T \quad \cdots \quad \mathbf{q}_{140,1}^T]$, $\mathbf{q}_{12}^T = [\mathbf{q}_{11,2}^T \quad \mathbf{q}_{12,2}^T \quad \cdots \quad \mathbf{q}_{140,2}^T]$, <br> in which <br> $\mathbf{q}_{10,1}^T = [9.313798\text{E-}19 \quad -2.0 \quad 12.0 \quad 1.672389\text{E-}17 \quad 0.0 \quad 0.0]$, <br> $\mathbf{q}_{11,1}^T = [-1.117232\text{E-}18 \quad -1.572635\text{E-}16 \quad 3.306817\text{E-}16 \quad 8.645323\text{E-}18 \quad 0.0 \quad 0.0]$, <br> ⋯⋯ |

$\mathbf{q}_{140,1}^{T} = [-4.520731\text{E-}19 \quad 2.358140\text{E-}18 \quad 5.101171\text{E-}16 \quad 3.161865\text{E-}18 \quad 0.0 \quad 0.0]$,

$\mathbf{q}_{11,2}^{T} = [1.171811\text{E-}18 \quad 3.756985\text{E-}16 \quad -7.473177\text{E-}16 \quad 1.534161\text{E-}17 \quad 0.0 \quad 0.0]$,

……

$\mathbf{q}_{140,2}^{T} = [1.292270\text{E-}18 \quad -2.366920\text{E-}16 \quad 3.668047\text{E-}16 \quad -1.166864\text{E-}18 \quad 0.0 \quad 0.0]$

c. $\mathbf{q}_{2}^{T}$

$\mathbf{q}_{2}^{T} = [\mathbf{q}_{21}^{T} \quad \mathbf{q}_{22}^{T}]$,

$\mathbf{q}_{21}^{T} = [\mathbf{q}_{20,1}^{T} \quad \mathbf{q}_{21,1}^{T} \quad \cdots \quad \mathbf{q}_{240,1}^{T}]$, $\mathbf{q}_{22}^{T} = [\mathbf{q}_{21,2}^{T} \quad \mathbf{q}_{22,2}^{T} \quad \cdots \quad \mathbf{q}_{240,2}^{T}]$,

in which

$\mathbf{q}_{20,1}^{T} = [1.423594\text{E-}18 \quad -2.0 \quad 12.0 \quad 1.672389\text{E-}17 \quad 0.0 \quad 0.0]$,

$\mathbf{q}_{21,1}^{T} = [1.371922\text{E-}18 \quad -2.072452\text{E-}16 \quad 5.442695\text{E-}17 \quad 8.645323\text{E-}18 \quad 0.0 \quad 0.0]$,

……

$\mathbf{q}_{240,1}^{T} = [1.674491\text{E-}18 \quad -1.761829\text{E-}18 \quad 4.322714\text{E-}16 \quad 3.161865\text{E-}18 \quad 0.0 \quad 0.0]$,

$\mathbf{q}_{21,2}^{T} = [-4.583626\text{E-}19 \quad 3.457093\text{E-}17 \quad -7.460165\text{E-}16 \quad 1.534161\text{E-}17 \quad 0.0 \quad 0.0]$,

……

$\mathbf{q}_{240,2}^{T} = [-1.736237\text{E-}18 \quad -2.090599\text{E-}16 \quad 4.645624\text{E-}16 \quad -1.166864\text{E-}18 \quad 0.0 \quad 0.0]$

d. $\mathbf{q}_{3}^{T}$

$\mathbf{q}_{3}^{T} = [0.0 \quad 0.0 \quad 0.0 \quad 0.0 \quad 4.0 \quad 0.0 \quad 0.0 \quad -24.0 \quad -24.0 \quad 0.0$
$\quad 0.0 \quad 0.0 \quad 144.0 \quad 0.0 \quad 0.0]$

7   a. $\mathbf{q}_{04}^{T}$

$\mathbf{q}_{04}^{T} = [1.339853\text{E-}17 \quad -5.925981\text{E-}18 \quad \cdots \quad 5.802761\text{E-}20 \quad -1.933270\text{E-}19]$

b. $\mathbf{q}_{1}^{T}$

$\mathbf{q}_{1}^{T} = [\mathbf{q}_{11}^{T} \quad \mathbf{q}_{12}^{T} \quad \cdots \quad \mathbf{q}_{140}^{T}]$,

in which

$\mathbf{q}_{11}^{T} = [-3.873731\text{E-}17 \quad -1.770942\text{E-}16 \quad 0.0 \quad -1.993994\text{E-}17 \quad 3.987988\text{E-}17 \quad 0.0]$,

$\mathbf{q}_{12}^{T} = [1.097717\text{E-}17 \quad 3.292473\text{E-}17 \quad 0.0 \quad 6.481345\text{E-}18 \quad -1.296269\text{E-}17 \quad 0.0]$,

……

$\mathbf{q}_{140}^{T} = [6.876119\text{E-}19 \quad 4.300453\text{E-}18 \quad 0.0 \quad 2.437254\text{E-}19 \quad -4.874507\text{E-}19 \quad 0.0]$,

c. $\mathbf{q}_{2}^{T}$

$\mathbf{q}_{2}^{T} = [\mathbf{q}_{21}^{T} \quad \mathbf{q}_{22}^{T} \quad \cdots \quad \mathbf{q}_{240}^{T}]$,

in which

$\mathbf{q}_{21}^{T} = [-3.706127\text{E-}17 \quad -1.760560\text{E-}16 \quad 0.0 \quad -1.993994\text{E-}17 \quad 3.987988\text{E-}17 \quad 0.0]$,

$\mathbf{q}_{22}^{T} = [9.086650\text{E-}18 \quad 2.494088\text{E-}17 \quad 0.0 \quad 6.481345\text{E-}18 \quad -1.296269\text{E-}17 \quad 0.0]$,

……

$\mathbf{q}_{240}^{T} = [5.016186\text{E-}19 \quad 6.721300\text{E-}18 \quad 0.0 \quad 2.437254\text{E-}19 \quad -4.874507\text{E-}19 \quad 0.0]$,

d. $\mathbf{q}_{3}^{T}$

$\mathbf{q}_{3}^{T} = [0.0 \quad 0.0 \quad 0.0 \quad 0.0 \quad 25.0 \quad 0.0 \quad 0.0 \quad 2.5 \quad 0.0 \quad 0.0 \quad -5.0 \quad 0.0$
$\quad 0.0 \quad 0.0 \quad 2.5 \quad 0.0 \quad -5.0 \quad 0.0 \quad 0.25 \quad -0.5 \quad -0.5 \quad 0.0 \quad 1.0 \quad 0.0]$

2. Errors in simultaneous approximation

Next, consider the sample functions in the form of trigonometrical functions. In Tables 6 and 7, we present the errors in the approximations respectively with the algebraical polynomial based composite Fourier series method (AP-CFSM) and the Fourier series direct-expansion method (FS-DEM). For the components up to fourth order in the (partial) derivative vectors of the sample functions, the overall errors of the algebraical polynomial based composite Fourier series method are 4-5 orders of magnitude smaller than those of the Fourier series direct-expansion method, the internal errors of the algebraical polynomial based composite Fourier series method are 3-5 orders of magnitude smaller than those of the Fourier series direct-expansion method, the boundary errors of the algebraical polynomial based composite Fourier series method are 4-7 orders of magnitude smaller than those of the Fourier series direct-expansion method, and the corner errors of the algebraical polynomial based composite Fourier series method are 2-5 orders of magnitude smaller than those of the Fourier series direct-expansion method. For the fifth order component(s) in the (partial) derivative vectors of the sample functions, the overall errors, internal errors and boundary errors of the algebraical polynomial based composite Fourier series method are respectively almost the same as or 1-3 orders of magnitude smaller than those of the Fourier series direct-expansion method, and the corner errors of the algebraical polynomial based composite Fourier series method are almost the same as or 2 orders of magnitude smaller than those of the Fourier series direct-expansion method. For the sixth order component(s) in the (partial) derivative vectors of the sample functions, the overall errors, internal errors and boundary errors of the algebraical polynomial based composite Fourier series method are respectively almost the same as those of the Fourier series direct-expansion method, and usually the corner errors of the algebraical polynomial based composite Fourier series method are almost the same as those of the Fourier series direct-expansion method (except for sample function 6). Therefore, it is verified that the algebraical polynomial based composite Fourier series method, as an improvement of the Fourier series direct-expansion method, is feasible for functions with varied boundary conditions.

Table 6: Errors in simultaneous composite Fourier series approximation of derivative vectors of one-dimensional sample functions.

| Approximation errors | Sample function 2 | | Sample function 4 | |
|---|---|---|---|---|
| | AP-CFSM | FS-DEM | AP-CFSM | FS-DEM |
| $\|e\|^4(u_{40})$ | 4.72516E-07 | 1.24047E-02 | 4.72516E-07 | 1.23744E-02 |
| $\|e_I\|^4(u_{40})$ | 4.72595E-07 | 1.22866E-02 | 4.72603E-07 | 1.23165E-02 |
| $\|e_B\|^4(u_{40})$ | 7.98557E-08 | 6.03144E-01 | 4.14033E-08 | 3.01591E-01 |
| $e^{(5)}(u_{40})$ | 2.06660E-04 | 2.06662E-04 | 2.06266E-04 | 2.06188E-04 |
| $e_I^{(5)}(u_{40})$ | 2.05129E-04 | 2.05131E-04 | 2.05511E-04 | 2.05434E-04 |
| $e_B^{(5)}(u_{40})$ | 7.85952E-03 | 7.85949E-03 | 3.97827E-03 | 3.97717E-03 |
| $e^{(6)}(u_{40})$ | 2.05366E-02 | 2.05368E-02 | 2.04864E-02 | 2.04799E-02 |
| $e_I^{(6)}(u_{40})$ | 2.03407E-02 | 2.03409E-02 | 2.03904E-02 | 2.03840E-02 |
| $e_B^{(6)}(u_{40})$ | 1.00000E+00 | 1.00000E+00 | 5.00000E-01 | 5.00000E-01 |

Table 7: Errors in simultaneous composite Fourier series approximation of partial derivative vectors of two-dimensional sample functions.

| Approxima-tion errors | Sample function 6 | | Sample function 8 | |
|---|---|---|---|---|
| | AP-CFSM | FS-DEM | AP-CFSM | FS-DEM |
| $\lVert e \rVert^4 (u_{40,40})$ | 5.18432E-06 | 2.01323E-02 | 2.15001E-07 | 1.80706E-02 |
| $\lVert e_I \rVert^4 (u_{40,40})$ | 2.43680E-06 | 7.86272E-03 | 2.10033E-07 | 1.08846E-02 |
| $\lVert e_B \rVert^4 (u_{40,40})$ | 4.67352E-05 | 3.21987E-01 | 3.25495E-07 | 1.94078E-01 |
| $\lVert e_C \rVert^4 (u_{40,40})$ | 2.62376E-03 | 2.00012E-01 | 1.45041E-06 | 2.00799E-01 |
| $\lvert e \rvert^5 (u_{40,40})$ | 1.30588E-03 | 1.88203E-02 | 6.18574E-05 | 1.66987E-02 |
| $\lvert e_I \rvert^5 (u_{40,40})$ | 5.83649E-04 | 6.49644E-03 | 4.39762E-05 | 1.06156E-02 |
| $\lvert e_B \rvert^5 (u_{40,40})$ | 1.04781E-02 | 3.18975E-01 | 4.84950E-04 | 1.64899E-01 |
| $\lvert e_C \rvert^5 (u_{40,40})$ | 8.62915E-01 | 5.00031E-01 | 1.98913E-03 | 2.50000E-01 |
| $\lvert e \rvert^6 (u_{40,40})$ | 9.95679E-02 | 2.10069E-02 | 5.36005E-03 | 1.63157E-02 |
| $\lvert e_I \rvert^6 (u_{40,40})$ | 4.75365E-02 | 8.77358E-03 | 3.41443E-03 | 9.51247E-03 |
| $\lvert e_B \rvert^6 (u_{40,40})$ | 9.41960E-01 | 3.23994E-01 | 5.15023E-02 | 1.83399E-01 |
| $\lvert e_C \rvert^6 (u_{40,40})$ | 4.41929E+01 | 3.39905E-17 | 2.04528E-01 | 1.44568E-01 |

## 6. Conclusions

Simultaneous approximation of functions and their (partial) derivatives is the motivation and challenge for the further development of new type methods for differential equations. In this paper, we perform thorough investigation on the problem of simultaneous approximation of the functions and their (partial) derivatives with composite Fourier series. The approximating function is invariantly sought as a linear combination of the corner function (only for the two-dimensional function), the boundary function(s) and the internal function. The corner function and the boundary function(s) here are introduced to take care of all the possible discontinuities associated with the original function and its relevant (partial) derivatives on the boundaries. Accordingly, the obtained composite Fourier series approximation uniformly converges to the original function. Also, this composite Fourier series is able to be simply differentiated, term by term, to yield uniformly convergent series expansions for the original function's (partial) derivatives up to 2$r$th order.

The current method can be presented in different versions corresponding to different supplementary functions involved in the corner function and the boundary function(s). Specifically, supplementary functions in the form of algebraical polynomials lead to the algebraical polynomial based composite Fourier series method. It is shown through numerical examples that this proposed method is universally applicable to the one-dimensional or two-dimensional functions with different periodical extensions and different boundary conditions. The reproducing property of complete algebraical polynomials of this method is also demonstrated theoretically and numerically.

The composite Fourier series method actually provides a practical way to approximate the unknown field functions and their (partial) derivatives simultaneously within the solution domain, which appears to have considerable potential for the analytical analysis of multiscale problems in science and engineering.

# Appendix A. Proof of the proposition that the algebraical polynomial based composite Fourier series method can reproduce a complete algebraical polynomial of degree $2r$

In this appendix, we summarize the third case in section 4 by the following Theorem A.1 and prove it in detail.

**Theorem A.1**

Let $r$ be a positive integer, then with Fourier series expansion of the two-dimensional function specified as the full-range composite Fourier series over the domain $[-a,a]\times[-b,b]$, the algebraical polynomial based composite Fourier series method can reproduce a two-dimensional complete algebraical polynomial of degree $2r$.

Proof. We define the set

$$\Omega = \{(j,l)| j,l = 0,1,2,\cdots, j+l \leq 2r\}, \tag{A.1}$$

and the algebraical polynomial

$$p_{jl}(x_1, x_2) = (x_1/a)^j (x_2/b)^l, \quad j,l = 0,1,2,\cdots, \tag{A.2}$$

then the set of the two-dimensional complete algebraical polynomials of degree $2r$ is given by

$$\mathbf{P}^{2r} = \{p_{jl}(x_1, x_2)\}_{(j,l)\in\Omega}. \tag{A.3}$$

For the algebraical polynomial based composite Fourier series method, we now consider the basis functions in the form of algebraical polynomials.

1. We define the sets

$$\Omega_{31} = \{(2j+2, 2l+2)| j,l = 0,1,2,\cdots, j+l \leq r-2\}, \tag{A.4}$$

$$\Omega_{32} = \{(2j+1, 2l+2)| j,l = 0,1,2,\cdots, j+l \leq r-2\}, \tag{A.5}$$

$$\Omega_{33} = \{(2j+2, 2l+1)| j,l = 0,1,2,\cdots, j+l \leq r-2\}, \tag{A.6}$$

$$\Omega_{34} = \{(2j+1, 2l+1)| j,l = 0,1,2,\cdots, j+l \leq r-1\}, \tag{A.7}$$

$$\Omega_3 = \Omega_{31} \cup \Omega_{32} \cup \Omega_{33} \cup \Omega_{34}, \tag{A.8}$$

then the set of the algebraical polynomials involved in the vector of basis functions, $\mathbf{\Phi}_3^T(x_1, x_2)$, is given by

$$\mathbf{P}_3 = \{p_{jl}(x_1, x_2)\}_{(j,l)\in\Omega_3}. \tag{A.9}$$

2. We define the set

$$\Omega_1 = \{(j,0)| j = 1,2,\cdots,2r\}, \tag{A.10}$$

then the set of the algebraical polynomials involved in the vector of basis functions, $\mathbf{\Phi}_1^T(x_1, x_2)$, is given by

$$\mathbf{P}_1 = \{p_{jl}(x_1, x_2)\}_{(j,l)\in\Omega_1}. \tag{A.11}$$

3. We define the set

$$\Omega_2 = \{(0,l)| l = 1,2,\cdots,2r\}, \tag{A.12}$$

then the set of the algebraical polynomials involved in the vector of basis functions, $\mathbf{\Phi}_2^T(x_1, x_2)$, is given by

$$\mathbf{P}_2 = \{p_{jl}(x_1, x_2)\}_{(j,l) \in \Omega_2}. \tag{A.13}$$

4. We define the set
$$\Omega_0 = \{(0,0)\}, \tag{A.14}$$
then the set of the algebraical polynomials involved in the vector of basis functions, $\mathbf{\Phi}_0^T(x_1, x_2)$, is given by
$$\mathbf{P}_0 = \{p_{jl}(x_1, x_2)\}_{(j,l) \in \Omega_0}. \tag{A.15}$$

Hence, the set of the algebraical polynomials reproduced by the algebraical polynomial based composite Fourier series method is given by
$$\mathbf{P}_0 \cup \mathbf{P}_1 \cup \mathbf{P}_2 \cup \mathbf{P}_3 = \{p_{jl}(x_1, x_2)\}_{(j,l) \in \Omega_0 \cup \Omega_1 \cup \Omega_2 \cup \Omega_3}. \tag{A.16}$$

Now we prove that
$$\Omega = \Omega_0 \cup \Omega_1 \cup \Omega_2 \cup \Omega_3. \tag{A.17}$$

1. Let $(j,l) \in \Omega_0 \cup \Omega_1 \cup \Omega_2 \cup \Omega_3$, then it is obvious that
$$(j,l) \in \Omega, \tag{A.18}$$
and
$$\Omega_0 \cup \Omega_1 \cup \Omega_2 \cup \Omega_3 \subset \Omega. \tag{A.19}$$

2. Let $(j,l) \in \Omega$.

① If $j=0$ and $l=0$, then
$$(j,l) \in \Omega_0 \subset \Omega_0 \cup \Omega_1 \cup \Omega_2 \cup \Omega_3. \tag{A.20}$$

② If $j=0$ and $l>0$, and since
$$j+l \leq 2r, \tag{A.21}$$
we have
$$1 \leq l \leq 2r, \tag{A.22}$$
hence
$$(j,l) \in \Omega_2 \subset \Omega_0 \cup \Omega_1 \cup \Omega_2 \cup \Omega_3. \tag{A.23}$$

③ If $j>0$ and $l=0$, and since
$$j+l \leq 2r, \tag{A.24}$$
we have
$$1 \leq j \leq 2r, \tag{A.25}$$
hence
$$(j,l) \in \Omega_1 \subset \Omega_0 \cup \Omega_1 \cup \Omega_2 \cup \Omega_3. \tag{A.26}$$

④ If $j>0$, $l>0$:

a. if $j$ and $l$ are even, then $j-2$ and $l-2$ are both even, and
$$j-2 \geq 0, \ l-2 \geq 0. \tag{A.27}$$

Let
$$j-2 = 2j_0 \text{ and } l-2 = 2l_0, \tag{A.28}$$
we have
$$j = 2j_0 + 2, \ l = 2l_0 + 2, \tag{A.29}$$
and
$$j_0 \geq 0, \ l_0 \geq 0, \tag{A.30}$$
and since $j+l \leq 2r$, we have
$$j_0 + l_0 \leq r-2, \tag{A.31}$$
hence

$$(j,l) \in \Omega_{31} \subset \Omega_3 \subset \Omega_0 \cup \Omega_1 \cup \Omega_2 \cup \Omega_3. \tag{A.32}$$

b. If $j$ is odd and $l$ is even, then $j-1$ and $l-2$ are both even, and
$$j-1 \geq 0, \quad l-2 \geq 0. \tag{A.33}$$

Let
$$j-1 = 2j_0, \quad l-2 = 2l_0, \tag{A.34}$$
then
$$j = 2j_0 + 1, \quad l = 2l_0 + 2, \tag{A.35}$$
and
$$j_0 \geq 0, \quad l_0 \geq 0, \tag{A.36}$$
and since $j+l \leq 2r$, we have
$$j_0 + l_0 \leq r - 2, \tag{A.37}$$
hence
$$(j,l) \in \Omega_{32} \subset \Omega_3 \subset \Omega_0 \cup \Omega_1 \cup \Omega_2 \cup \Omega_3. \tag{A.38}$$

c. If $j$ is even and $l$ is odd, then $j-2$ and $l-1$ are both even, and
$$j-2 \geq 0, \quad l-1 \geq 0. \tag{A.39}$$

Let
$$j-2 = 2j_0, \quad l-1 = 2l_0, \tag{A.40}$$
then
$$j = 2j_0 + 2, \quad l = 2l_0 + 1, \tag{A.41}$$
and
$$j_0 \geq 0, \quad l_0 \geq 0, \tag{A.42}$$
and since $j+l \leq 2r$, we have
$$j_0 + l_0 \leq r - 2, \tag{A.43}$$
hence
$$(j,l) \in \Omega_{33} \subset \Omega_3 \subset \Omega_0 \cup \Omega_1 \cup \Omega_2 \cup \Omega_3. \tag{A.44}$$

d. If $j$ and $l$ are odd, then $j-1$ and $l-1$ are both even, and
$$j-1 \geq 0, \quad l-1 \geq 0. \tag{A.45}$$

Let
$$j-1 = 2j_0, \quad l-1 = 2l_0, \tag{A.46}$$
then
$$j = 2j_0 + 1, \quad l = 2l_0 + 1, \tag{A.47}$$
and
$$j_0 \geq 0, \quad l_0 \geq 0, \tag{A.48}$$
and since $j+l \leq 2r$, we have
$$j_0 + l_0 \leq r - 1, \tag{A.49}$$
hence
$$(j,l) \in \Omega_{34} \subset \Omega_3 \subset \Omega_0 \cup \Omega_1 \cup \Omega_2 \cup \Omega_3. \tag{A.50}$$

By combining Eqs. (A.20), (A.23), (A.26), (A.32), (A.38), (A.44) and (A.50), we have
$$\Omega \subset \Omega_0 \cup \Omega_1 \cup \Omega_2 \cup \Omega_3. \tag{A.51}$$

According to Eq. (A.19) and Eq. (A.51), we know that Eq. (A.17) is true and accordingly
$$\mathbf{P}_0 \cup \mathbf{P}_1 \cup \mathbf{P}_2 \cup \mathbf{P}_3 = \mathbf{P}^{2r}. \tag{A.52}$$

Therefore, the theorem follows. □